# Cruise Controllers for Lane-Free Ring-Roads based on Control Lyapunov Functions

**Dionysis Theodosis[*], Iasson Karafyllis[***], and Markos Papageorgiou[*],[**]**

[*] Dynamic Systems and Simulation Laboratory,
Technical University of Crete, Chania, 73100, Greece
[**] Faculty of Maritime and Transportation, Ningbo University,
Ningbo, China
emails: dtheodosis@dssl.tuc.gr , markos@dssl.tuc.gr
[***] Dept. of Mathematics, National Technical University of Athens,
Zografou Campus, 15780, Athens, Greece,
emails: iasonkar@central.ntua.gr , iasonkaraf@gmail.com

**Abstract**

The paper introduces novel families of cruise controllers for autonomous vehicles on lane-free ring-roads. The design of the cruise controllers is based on the appropriate selection of a Control Lyapunov Function expressed on measures of the energy of the system with the kinetic energy expressed in ways similar to Newtonian or relativistic mechanics. The derived feedback laws (cruise controllers) are decentralized (per vehicle), as each vehicle determines its control input based on: (i) its own state; (ii) either only the distance from adjacent vehicles (inviscid cruise controllers) or the state of adjacent vehicles (viscous cruise controllers); and (iii) its distance from the boundaries of the ring-road. A detailed analysis of the differences and similarities between lane-free straight roads and lane-free ring-roads is also presented.

## 1. Introduction

Designing safe and efficient control strategies for autonomous vehicles constitutes a challenging topic, addressed for a variety of driving circumstances, including Adaptive Cruise Control (ACC) and Cooperative ACC (CACC) systems where vehicles can communicate with each other, see for instance [1], [8], [20], [23], [36] and references therein. The impact of autonomous vehicles on traffic flow is typically studied on a string of vehicles that operates on a single-lane road or on a ring-road. The objective is then to design feedback laws so that vehicles can adjust their speed to the speed of the leader while keeping a safe distance from the following vehicle, see [2], [5], [9], [13], [21], [25], [34].

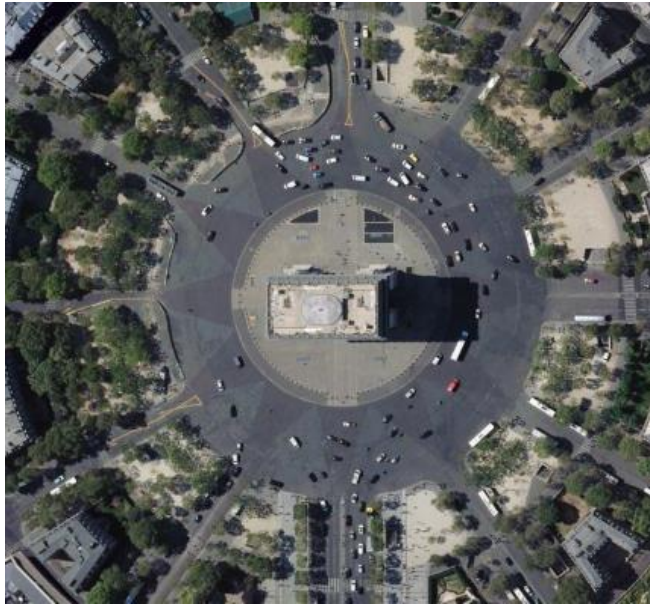

**Figure 1:** Place Charles-de-Gaulle (photo by National Institute of Geographic and Forest Information, source: Wikipedia)

Ring-roads are of particular interest since they may give rise to phantom traffic jams when the average density is higher than the critical density, see [3], [31], [35]. In this context, the ability of connected and automated vehicles on ring-roads to dissipate traffic waves have been intensively studied and reported in both microscopic and macroscopic traffic flow models, see [6], [7], [11], [12], [14], [15]. With the advancement of vehicular technology and the emergence of highly automated vehicles, new directions have been proposed, [24], where vehicles are not bound to traffic lanes but can freely move on the two-dimensional surface of the road and can influence the movement of adjacent vehicles all around them through sensors and communication, see [17], [18], [22], [37].

While cruise controllers for autonomous vehicles on ring-roads and roundabouts have a long history (see for instance [5], [9], [13], [21], [25], [26], [28], [30], [33], [38]), they are based on a single or double-lane road and are therefore not suitable to capture the complexity of lane free ring-roads, such as the famous Charles-de-Gaulle Place roundabout in Paris, France, (Figure 1), see also [22] for a vehicle movement strategy on the latter roundabout. In this paper, we extend the Control Lyapunov Function (CLF) methodologies presented in [18] to derive families of cruise controllers for autonomous vehicles that operate on the two-dimensional surface of lane-free ring-roads. The CLFs are based on measures of the total energy of the system and act as size functions (Lemma 1 and Lemma 2) guaranteeing that the closed-loop system is well-posed. By expressing the kinetic energy in ways similar to Newtonian or relativistic mechanics, two respective families of cruise controllers are obtained that satisfy the following properties globally (Theorem 1 and Theorem 2): (i) there are no collisions among vehicles or with the inner and outer boundaries of the ring-road; (ii) the speeds of all vehicles are always positive and remain below a given speed limit; (iii) the angular speeds of all vehicles converge to a given angular speed set-point; and the accelerations, the relative orientations (the deviation of the heading angle from the tangent of the circle), and rates of change of the relative orientations of all vehicles tend to zero. The proposed families of cruise controllers are decentralized (per vehicle), and each vehicle only has access to its own state, the distance from the boundaries of the ring-road and either the relative positions from adjacent vehicles (inviscid cruise controllers) or the relative speeds and relative positions from adjacent vehicles (viscous cruise controllers).

By designing cruise controllers for both the case of a ring-road and the case of a straight road (see [18]), it is possible to consider any road represented by a closed or an open non-self-intersecting curve by selecting a suitable change of coordinates to:

- Transform any curved road of constant width and infinite length into a straight road of constant width;
- Transform any road of constant width represented by a closed curve into a ring-road of constant width.

It should be highlighted that, since the straight road and the ring-road are "topologically" different, it is not possible to find a change of coordinates to transform a straight road into a ring-road, and therefore the cruise controllers proposed in this paper are not a direct consequence of the cruise controllers of the straight road case.

The structure of the paper is as follows. Section 2 is devoted to the description of the employed vehicle model for the ring-road. Section 3 presents the families of the proposed cruise controllers and the statements of our main results. In Section 4, the differences and similarities between the cruise controllers for the ring-road and the straight road are addressed. Simulation examples are included in Section 5 that demonstrate the properties of the families of cruise controllers. All proofs of the main results are provided in Section 6, while some concluding remarks are given in Section 7.

**Notation.** Throughout this paper, we adopt the following notation.

∗ $\Re_+ := [0, +\infty)$ denotes the set of non-negative real numbers.

∗ By $|x|$ we denote both the Euclidean norm of a vector $x \in \Re^n$ and the absolute value of a scalar $x \in \Re$. By $x'$ we denote the transpose of a vector $x \in \Re^n$. By $|x|_\infty = \max\{|x_i|, i = 1, ..., n\}$ we denote the infinity norm of a vector $x = (x_1, x_2, ..., x_n)' \in \Re^n$.

∗ Let $A \subseteq \Re^n$ be an open set. By $C^0(A, \Omega)$, we denote the class of continuous functions on $A \subseteq \Re^n$, which take values in $\Omega \subseteq \Re^m$. By $C^k(A; \Omega)$, where $k \geq 1$ is an integer, we denote the class of functions on $A \subseteq \Re^n$ with continuous derivatives of order $k$, which take values in $\Omega \subseteq \Re^m$. When $\Omega = \Re$ the we write $C^0(A)$ or $C^k(A)$. For a function $V \in C^1(A; \Re)$, the gradient of $V$ at $x \in A \subseteq \Re^n$, denoted by $\nabla V(x)$, is the row vector $\left[ \frac{\partial V}{\partial x_1}(x) \cdots \frac{\partial V}{\partial x_n}(x) \right]$.

∗ Given $f : S \to \Re$, $S \subset \Re$, we denote by $\arg\min_{x \in S}(f(x))$ the set of elements in $S$ that achieve the global minimum in $S$, i.e., $\arg\min_{x \in S} f(x) = \{x \in S : f(s) \geq f(x), \text{ for all } s \in S\}$.

∗ We denote by $dist(x, A)$ the Euclidean distance of the point $x \in \Re^n$ from the set $A \subset \Re^n$, i.e., $dist(x, A) = \inf\{|x - y| : y \in A\}$.

## 2. Model and Problem Formulation

The bicycle kinematic model ([27]) has been widely used to describe the motion of vehicles on a straight road

$$\begin{aligned}
\dot{x}_i &= v_i \cos(\theta_i) \\
\dot{y}_i &= v_i \sin(\theta_i) \\
\dot{\theta}_i &= \frac{v_i}{\sigma_i} \tan(\delta_i) \\
\dot{v}_i &= F_i
\end{aligned} \tag{2.1}$$

where $(x_i, y_i)$, $i=1,\ldots,n$ is the reference point of the $i$-th vehicle placed at the midpoint of the rear axle of the vehicle in an inertial frame with Cartesian coordinates $(X,Y)$, $v_i$ is the speed of the vehicle at the point $(x_i, y_i)$, $\theta_i$ is the heading angle (orientation) with respect to the X axis, $\sigma_i$ is the length of the vehicle $i$ (a constant), $F_i$ is the acceleration and $\delta_i$ is the steering angle of the front wheels relative to the orientation $\theta_i$ of the vehicle.

To describe the motion of a vehicle on a ring-road of inner radius $R_{in} > 0$ and outer radius $R_{out} > R_{in}$ centered at $(0,0)$, we let $r_i \in (R_{in}, R_{out})$ be the distance of the reference point of vehicle $i$ from $(0,0)$, and $\varphi_i$ be the angular coordinate (the angle of the reference point of vehicle $i$ from the reference direction $X$). Then, by using the change of coordinates $x_i = r_i \cos(\varphi_i)$, $y_i = r_i \sin(\varphi_i)$ we obtain the bicycle model in polar coordinates, as appropriate for movement on a ring-road:

$$\begin{aligned}
\dot{r}_i &= -v_i \sin(s_i) \\
\dot{\varphi}_i &= \frac{v_i}{r_i} \cos(s_i) \\
\dot{s}_i &= \frac{v_i}{\sigma_i} \tan(\delta_i) - \frac{v_i}{r_i} \cos(s_i) \\
\dot{v}_i &= F_i
\end{aligned} \tag{2.2}$$

for $i = 1,\ldots,n$, where $v_i \in (0, v_{\max})$ is the speed of the $i$-th vehicle at the point $(r_i, \varphi_i)$, $v_{\max} > 0$ denotes the road speed limit, $s_i = \theta_i - \varphi_i - \frac{\pi}{2}$ is the relative orientation, i.e., the deviation of the heading angle from the tangent of the circle with radius $r_i$ and center $(0,0)$.

Let $\omega^* \in \left(0, \frac{v_{\max}}{R_{out}}\right)$ be given (the angular speed set point) and define the set

$$S := (R_{in}, R_{out})^n \times \mathfrak{R}^n \times (-\Theta, \Theta)^n \times (0, v_{\max})^n \tag{2.3}$$

where $\Theta \in \left(0, \frac{\pi}{2}\right)$ is a given angle (maximum deviation of heading angle from the tangent of the ring road) that satisfies

$$\cos(\Theta) > \frac{R_{out}\omega^*}{v_{max}}. \tag{2.4}$$

We define the distance between vehicles by

$$d_{i,j} := \sqrt{p_{i,j}(r_i - r_j)^2 + 2r_i r_j (1 - \cos(\varphi_i - \varphi_j))}, \text{ for } i, j = 1,...,n \tag{2.5}$$

where $p_{i,j} > 0$ are weight parameters that satisfy $p_{i,j} = p_{j,i}$ for all $i, j = 1,...,n$. Notice that when $p_{i,j} = 1$ the distance metric defined by (2.5) coincides with the usual Euclidean distance metric (i.e., $d_{i,j}$ is equal to the usual Euclidean distance between the points $(r_i, \varphi_i)$ and $(r_j, \varphi_j)$ in polar coordinates). Let

$$w = (r_1,...,r_n, \varphi_1,...,\varphi_n, s_1,...,s_n, v_1,...,v_n)' \in \Re^{4n} \tag{2.6}$$

Due to the various constraints presented above, the state space of the model (2.2) is

$$\Omega := \{ w \in S : d_{i,j} > L_{i,j}, i, j = 1,...,n, j \neq i \} \tag{2.7}$$

Where $L_{i,j}$, $i, j = 1,...,n$, $i \neq j$, are positive constants and represent the minimum distance between a vehicle $i$ and a vehicle $j$, with $L_{i,j} = L_{j,i}$ for $i, j = 1,...,n$, $i \neq j$. Notice that the state-space $\Omega$ is not a linear subspace of $\Re^{4n}$ but an open set.

**Problem Statement:** Design cruise controllers for vehicles operating on lane-free ring-roads that satisfy the following properties:

(P1) Well-posedness requirement: For each $w(0) \in \Omega$, there exists a unique solution $w(t) \in \Omega$ defined for all $t \geq 0$. According to (2.7), this requirement implies that there are no collisions between vehicles (since $d_{i,j}(t) > L_{i,j}$ for $t \geq 0$, $i, j = 1,...,n$, $j \neq i$) or with the boundary of the road (since $r_i(t) \in (R_{in}, R_{out})$ for $t \geq 0$); the speeds of all vehicles are always positive and remain below the given speed limit (since $v_i(t) \in (0, v_{max})$ for all $t \geq 0$); and the orientation of each vehicle is always bounded by the given value $\Theta \in \left(0, \frac{\pi}{2}\right)$ (since $s_i(t) \in (-\Theta, \Theta)$ for $t \geq 0$).

(P2) Asymptotic requirement: The relative orientation of each vehicle satisfies $\lim_{t \to +\infty}(s_i(t)) = 0$ for $i = 1,...,n$, and the angular speed of all vehicles satisfy $\lim_{t \to +\infty}\left(\frac{v_i(t)}{r_i(t)}\right) = \omega^*$, $i = 1,...,n$, for a given angular speed set-point $\omega^* \in \left(0, \frac{v_{max}}{R_{out}}\right)$. Moreover, the accelerations and the rate of change of the relative orientations of all vehicles tend to zero, i.e., $\lim_{t \to +\infty}(F_i(t)) = 0$, and $\lim_{t \to +\infty}(\dot{s}_i(t)) = 0$, for $i = 1,...,n$.

# 3. Main Results

## 3.1 Preliminaries

In this section we present two families of cruise controllers for vehicles operating on lane-free ring-roads that address properties (P1) and (P2).

Let $V_{i,j} : (L_{i,j}, +\infty) \to \Re_+$, $U_i : (R_{in}, R_{out}) \to \Re_+$, $i, j = 1, ..., n$, $j \neq i$ be $C^2$ functions and $\kappa_{i,j} : (L_{i,j}, +\infty) \to \Re_+$ $i, j = 1, ..., n$, $j \neq i$ be $C^1$ functions that satisfy the following properties

$$\lim_{d \to L_{i,j}^+} \left( V_{i,j}(d) \right) = +\infty, \ i, j = 1, ..., n, j \neq i \tag{3.1}$$

$$V_{i,j}(d) = 0, \text{ for all } d \geq \lambda, \ i, j = 1, ..., n, j \neq i \tag{3.2}$$

$$V_{i,j}(d) \equiv V_{j,i}(d), \ i, j = 1, ..., n, j \neq i \tag{3.3}$$

$$\lim_{r \to R_{in}^+} \left( U_i(r) \right) = +\infty, \ \lim_{r \to R_{out}^-} \left( U_i(r) \right) = +\infty, \ i = 1, ..., n \tag{3.4}$$

$$\kappa_{i,j}(d) \equiv \kappa_{j,i}(d), \ i, j = 1, ..., n, j \neq i \tag{3.5}$$

$$\kappa_{i,j}(d) = 0, \text{ for all } d \geq \lambda, \ i, j = 1, ..., n, j \neq i \tag{3.6}$$

where $\lambda$ is a positive constant that satisfies

$$\lambda > \max \{ L_{i,j}, i, j = 1, ..., n, i \neq j \} \tag{3.7}$$

The functions $V_{i,j}$ and $U_i$ are potential functions which are used to avoid collisions between vehicles and road boundary violation, respectively (see [18], [36]). Condition (3.3) implies that if a vehicle $i$ exerts a force to vehicle $j$, then vehicle $j$ exerts the opposite force to vehicle $i$. In addition, properties (3.2) and (3.6) will allow to design decentralized controllers, which use real-time information (such as relative positions, relative speeds, and relative orientations) only from adjacent vehicles that are located at a distance less than $\lambda > 0$. Finally, the functions $\kappa_{i,j}$ satisfying (3.5), (3.6) will be used to introduce a viscous-like behavior of the vehicles. Notice that $\kappa_{i,j}$ are not necessarily potential functions as it is not assumed that they satisfy the limit property $\lim_{d \to L_{i,j}^+} \left( \kappa_{i,j}(d) \right) = +\infty$.

In the following sections, we will use a Control Lyapunov Function (CLF) methodology and the potential functions above to derive families of cruise controllers for autonomous vehicles on lane-free ring-roads. The construction of the Lyapunov function is based on the measures of the total energy of the system. Depending on how the kinetic energy of the system is expressed, we obtain two different families of cruise controllers. If the kinetic energy is expressed in a fashion similar to that of Newtonian mechanics, we call the corresponding controller a *Newtonian Cruise Controller* (NCC); while, when the kinetic energy is expressed in terms similar to that of relativistic mechanics, we call the corresponding controller a *Pseudo-Relativistic Cruise Controller* (PRCC). The main difference of those two approaches is that, in relativistic mechanics, the kinetic energy increases to infinity when an object's speed approaches the speed of light, instead of which the maximum allowed speed $v_{max}$ is used here; while the kinetic energy in Newtonian mechanics continues to increase without bound as the speed of an object increases.

Finally, when $\kappa_{i,j}(d) \equiv 0$ for all $i, j = 1,...,n$., $i \neq j$, we call the controller *"inviscid"* since the corresponding macroscopic model does not contain a viscosity term; otherwise, the corresponding controller is called *"viscous"*, see [18] for the corresponding macroscopic models.

## 3.2 Newtonian Cruise Controller (NCC)

The CLF in this case is given by the formula

$$H(w) := \frac{1}{2}\sum_{i=1}^{n}\left(\frac{v_i}{r_i}\cos(s_i) - \omega^*\right)^2 + \frac{b}{2}\sum_{i=1}^{n}v_i^2\sin^2(s_i)$$
$$+ \sum_{i=1}^{n}U_i(r_i) + \frac{1}{2}\sum_{i=1}^{n}\sum_{j \neq i}V_{i,j}(d_{i,j}) + A\sum_{i=1}^{n}\left(\frac{1}{\cos(s_i) - \cos(\Theta)} - \frac{1}{1 - \cos(\Theta)}\right)$$
(3.8)

where $A > 0$, $b > \frac{1}{R_{in}^2}$ are parameters of the Lyapunov function. The function $H$ in (3.8) is based on the total mechanical energy of the system of $n$ vehicles. Specifically, the first two terms ($\frac{1}{2}\sum_{i=1}^{n}\left(\frac{v_i}{r_i}\cos(s_i) - \omega^*\right)^2 + \frac{b}{2}\sum_{i=1}^{n}v_i^2\sin^2(s_i)$) are related to the rotational kinetic energy of the system of $n$ vehicles relative to an observer moving along the ring-road with angular speed equal to $\omega^*$ (as in classical mechanics). The sum of the third and fourth terms ($\sum_{i=1}^{n}U_i(r_i) + \frac{1}{2}\sum_{i=1}^{n}\sum_{j \neq i}V_{i,j}(d_{i,j})$), is related to the potential energy of the system. Finally, the last term of (3.8) ($A\sum_{i=1}^{n}\left(\frac{1}{\cos(s_i) - \cos(\Theta)} - \frac{1}{1 - \cos(\Theta)}\right)$) is a penalty term that blows up when $s_i \to \pm\Theta$.

The following lemma shows that the CLF (3.8) has certain properties of a size function (see [29]), however, it is not a (global) size function.

**Lemma 1:** *Let constants $R_{out} > R_{in} > 0$, $A > 0$, $v_{max} > 0$, $\omega^* \in \left(0, \frac{v_{max}}{R_{out}}\right)$, $L_{i,j} > 0$, $i, j = 1,...,n$, $i \neq j$, $\lambda > 0$ that satisfies (3.7); $\Theta \in \left(0, \frac{\pi}{2}\right)$ that satisfies (2.4); and define the function $H : \Omega \to \Re_+$ by means of (3.8), where $\Omega$ is given by (2.7). Then, there exist constants $\xi_i \in (R_{in}, R_{out})$, non-decreasing functions $\eta_i : \Re_+ \to [\xi_i, R_{out})$, $\omega : \Re_+ \to [0, \Theta)$, $i = 1,...,n$, non-increasing functions $\zeta_i : \Re_+ \to (R_{in}, \xi_i]$, $i = 1,...,n$, and, for each pair $i, j = 1,...,n$, $i \neq j$, there exist non-increasing functions $\rho_{i,j} : \Re_+ \to (L_{i,j}, \lambda]$ with $\rho_{i,j}(s) \equiv \rho_{j,i}(s)$, such that the following implications hold:*

$$w \in \Omega \Rightarrow \zeta_i(H(w)) \leq r_i \leq \eta_i(H(w)), |s_i| \leq \omega(H(w)), d_{i,j} \geq \rho_{i,j}(H(w))$$

*for $i, j = 1,...,n, j \neq i$.* (3.9)

Let $g_k : \Re \to \Re$ $(k=1,2)$ be given $C^1$ non-decreasing functions and $f : \Re \to \Re_+$ be a $C^1$ function that satisfies

$$f(x) \geq \max(0, x), \text{ for all } x \in \Re. \tag{3.10}$$

Based on the CLF (3.8), the Newtonian Cruise Controller (NCC) for (2.2) is given by the following equations for all $w \in \Omega$:

$$F_i = -k_i(w)\left(v_i - \frac{r_i \omega^*}{\cos(s_i)}\right) - \frac{r_i \omega^*}{\cos(s_i)}\left(\Phi_i(w) - G_i(w)\right), \text{ for } i=1,\dots,n \tag{3.11}$$

$$\delta_i = \tan^{-1}\left(\frac{\sigma_i}{r_i}\cos(s_i) - \frac{\sigma_i}{v_i a(r_i, s_i, v_i)}\left(\mu_2 \sin(s_i) + (bF_i \sin(s_i) + \Lambda_i(w))v_i - M_i(w)\right)\right),$$
$$\text{for } i=1,\dots,n \tag{3.12}$$

where

$$k_i(w) = \mu_1 + \Phi_i(w) - G_i(w) + f\left(-\frac{v_{\max} \cos(s_i)}{v_{\max} \cos(s_i) - r_i \omega^*}\left(\Phi_i(w) - G_i(w)\right)\right),$$
$$\text{for } i=1,\dots,n \tag{3.13}$$

and $\mu_1, \mu_2 > 0$ are positive constants. Moreover, we have

$$a(r,s,v) := \left(b - \frac{1}{r^2}\right)v^2 \cos(s) + \omega^* \frac{v}{r} + \frac{A}{(\cos(s) - \cos(\Theta))^2}$$
$$\text{for } (r,s,v) \in (R_{in}, R_{out}) \times (-\Theta, \Theta) \times (0, v_{\max}) \tag{3.14}$$

$$\Lambda_i(w) := \left(\frac{v_i}{r_i}\cos(s_i) - \omega^*\right)\frac{v_i}{r_i^2}\cos(s_i) - U_i'(r_i)$$
$$- \sum_{j \neq i}\left(p_{i,j}(r_i - r_j) + r_j(1 - \cos(\varphi_i - \varphi_j))\right)\frac{V_{i,j}'(d_{i,j})}{d_{i,j}} \tag{3.15}$$

$$\Phi_i(w) := \frac{r_i}{\omega^*}\sum_{j \neq i} V_{i,j}'(d_{i,j})\frac{r_j \sin(\varphi_i - \varphi_j)}{d_{i,j}} \tag{3.16}$$

$$G_i(w) := \frac{1}{\omega^*}\sum_{j \neq i} \kappa_{i,j}(d_{i,j})\left(g_1\left(\frac{v_j}{r_j}\cos(s_j)\right) - g_1\left(\frac{v_i}{r_i}\cos(s_i)\right)\right) \tag{3.17}$$

and

$$M_i(w) := \sum_{j \neq i} \kappa_{i,j}(d_{i,j})\left(g_2(\sin(s_j)) - g_2(\sin(s_i))\right)$$
$$\text{for } w \in \Omega, \, i=1,\dots,n \tag{3.18}$$

It should be noticed that when the NCC is inviscid, i.e., when either the functions $\kappa_{i,j}$ are zero or the functions $g_k$ are constant, then, the only real-time measurement requirements of each vehicle are its own state and the distances of its adjacent vehicles (not their speeds, not their orientations). In contrast, when the NCC is viscous, i.e., when the functions $\kappa_{i,j}$ are not zero and the functions $g_k$ are not constant, then the real-time measurement requirements of each vehicle are its own state and the states of its adjacent vehicles (including their speeds and their orientations). The term $k_i(w)$ in the acceleration $F_i(t)$ given in (3.11) is a state-dependent controller gain, which guarantees that the speed of each vehicle will remain positive and less than the speed limit $v_{\max}$. Finally, notice that properties (3.2) and (3.6) guarantee that the feedback laws (3.11) and (3.12), are decentralized (per vehicle) and depend only on the relative positions from adjacent vehicles for the inviscid case, or on the states of adjacent vehicles for the viscous case.

## 3.3. Pseudo-Relativistic Cruise Controller (PRCC)

The CLF in this case is given by the formula

$$H_R(w) := \frac{1}{2}\sum_{i=1}^{n} \frac{\left(\frac{v_i}{r_i}\cos(s_i) - \omega^*\right)^2 + bv_i^2 \sin^2(s_i)}{(v_{\max} - v_i)v_i} \quad (3.19)$$

$$+ \sum_{i=1}^{n} U_i(r_i) + \frac{1}{2}\sum_{i=1}^{n}\sum_{j \neq i} V_{i,j}(d_{i,j}) + A\sum_{i=1}^{n}\left(\frac{1}{\cos(s_i) - \cos(\Theta)} - \frac{1}{1-\cos(\Theta)}\right)$$

where $A > 0$, $b > \frac{1}{R_{in}^2}$ are parameters of the Lyapunov function. Notice that the kinetic energy term in $H_R$ (i.e., the term $\frac{1}{2}\sum_{i=1}^{n} \frac{\left(\frac{v_i}{r_i}\cos(s_i) - \omega^*\right)^2 + bv_i^2 \sin^2(s_i)}{(v_{\max} - v_i)v_i}$) is similar to the kinetic energy of a system of $n$ particles in relativistic mechanics, with speed limits 0 and $v_{\max}$ in place of $-c$ and $c$, where $c$ is the speed of light, which are the speed limits in relativistic mechanics. In relativistic mechanics, the kinetic energy increases to infinity when the speed of an object approaches (in absolute value) the speed of light, which indicates that no object with mass can reach the speed of light. Analogously, in (3.19), the kinetic energy term grows to infinity as the speed of a vehicle approaches zero or the maximum speed $v_{\max}$, thus restricting the speed of vehicles in $(0, v_{\max})$. As in the case of (3.8), the sum of the terms $\sum_{i=1}^{n} U_i(r_i) + \frac{1}{2}\sum_{i=1}^{n}\sum_{j \neq i} V_{i,j}(d_{i,j})$ is related to the potential energy of the system; and the last term of (3.19) ( $A\sum_{i=1}^{n}\left(\frac{1}{\cos(s_i) - \cos(\Theta)} - \frac{1}{1-\cos(\Theta)}\right)$) is a penalty term that blows up when $s_i \to \pm\Theta$.

The following lemma shows that the CLF (3.19) is a size function (see [29]) for the state space $\Omega$ defined by (2.7).

**Lemma 2:** *Let constants $A > 0$, $v_{max} > 0$, $\omega^* \in \left(0, \dfrac{v_{max}}{R_{out}}\right)$, $L_{i,j} > 0$, $i, j = 1, \ldots, n$, $i \neq j$, $\lambda > 0$ that satisfies (3.7), $\Theta \in \left(0, \dfrac{\pi}{2}\right)$ that satisfies (2.5), and define the function $H_R : \Omega \to \Re_+$ by means of (3.19), where $\Omega$ is given by (2.7). Then, there exist constants $\bar{v} \in (0, v_{max})$, $\xi_i \in (R_{in}, R_{out})$, $i = 1, \ldots, n$, non-decreasing functions $\eta_i : \Re_+ \to [\xi_i, R_{out})$, $i = 1, \ldots, n$, $\psi : \Re_+ \to [\bar{v}, v_{max})$, $\omega : \Re_+ \to [0, \Theta)$, non-increasing functions $\zeta_i : \Re_+ \to (R_{in}, \xi_i]$, $i = 1, \ldots, n$, $\ell : \Re_+ \to (0, \bar{v}]$, and, for each pair $i, j = 1, \ldots, n$, $i \neq j$, there exist non-increasing functions $\rho_{i,j} : \Re_+ \to (L_{i,j}, \lambda]$ with $\rho_{i,j}(s) \equiv \rho_{j,i}(s)$, such that the following implications hold:*

$$w \in \Omega \Rightarrow \begin{array}{l} \ell(H_R(w)) \leq v_i \leq \psi(H_R(w)), |s_i| \leq \omega(H_R(w)), \\ \zeta_i(H_R(w)) \leq r_i \leq \eta_i(H_R(w)), d_{i,j} \geq \rho_{i,j}(H_R(w)) \end{array} \quad \text{for } i, j = 1, \ldots, n, j \neq i. \quad (3.20)$$

Let $g_k : \Re \to \Re$ ($k = 1, 2$) be given $C^1$ non-decreasing functions and let also $f_j : \Re \to \Re$, $j = 1, 2$, be $C^1$ functions that satisfy:

$$f_j(0) = 0 \text{ and } x f_j(x) > 0, \text{ for } x \neq 0, \; j = 1, 2 \quad (3.21)$$

The Pseudo-Relativistic Cruise Controller (PRCC) for (2.2) that correspond to the CLF (3.19) is given by the following equations for all $w \in \Omega$:

$$F_i = -\dfrac{1}{q(r_i, s_i, v_i)} \left( f_1\left( \dfrac{v_i}{r_i} \cos(s_i) - \omega^* \right) + \omega^* (\Phi_i(w) - G_i(w)) \right), \text{ for } i = 1, \ldots, n \quad (3.22)$$

$$\delta_i = \tan^{-1}\left( \dfrac{\sigma_i}{r_i} \cos(s_i) - \dfrac{\sigma_i}{\gamma(r_i, s_i, v_i) v_i} \left( f_2(\sin(s_i)) - (\zeta_i(s_i, v_i) F_i + Z_i(w)) v_i - M_i(w) \right) \right),$$
$$\text{for } i = 1, \ldots, n \quad (3.23)$$

where $\Phi_i$, $G_i$ and $M_i$ are given by (3.16), (3.17) and (3.18), respectively, and

$$q(r, s, v) = \dfrac{v_{max} v \cos(s) - 2rv\omega^* + r\omega^* v_{max}}{2r(v_{max} - v)^2 v^2} \quad (3.24)$$

$$\gamma(r, s, v) = \dfrac{A}{(\cos(s) - \cos(\Theta))^2} + \dfrac{v \cos(s)}{v_{max} - v}\left( b - \dfrac{1}{r^2} \right) + \dfrac{\omega^*}{r(v_{max} - v)} \quad (3.25)$$

$$\zeta(s, v) = \dfrac{b v_{max} \sin(s)}{2(v_{max} - v)^2 v} \quad (3.26)$$

for $(r, s, v) \in (R_{in}, R_{out}) \times (-\Theta, \Theta) \times (0, v_{max})$, $b > \dfrac{1}{R_{in}^2}$ and

$$Z_i(w) := \left(\frac{v_i}{r_i}\cos(s_i) - \omega^*\right)\frac{\cos(s_i)}{(v_{\max} - v_i)r_i^2} - U_i'(r_i)$$
$$-\sum_{j \neq i}\left(p_{i,j}(r_i - r_j) + r_j(1 - \cos(\varphi_i - \varphi_j))\right)\frac{V'_{i,j}(d_{i,j})}{d_{i,j}} \quad (3.27)$$

for $w \in \Omega$. Notice that the assumption $b > \frac{1}{R_{in}^2}$ guarantees that $\gamma(r,s,v) > 0$ for all $(r,s,v) \in (R_{in}, R_{out}) \times (-\Theta, \Theta) \times (0, v_{\max})$.

The pseudo-relativistic feedback laws (3.22) and (3.23) are derived by using the CLF (3.19) which is also a size function. Compared to the Newtonian controller (3.11), (3.12), the controller (3.21), (3.22) is simpler, since it does not use state-dependent controller gains to restrict the speed in $(0, v_{\max})$ (due to the properties of the size function $H_R$). Finally, notice that properties (3.2) and (3.6) guarantee that the feedback laws (3.22) and (3.23), are decentralized (per vehicle) and depend only on adjacent vehicles' states, namely vehicles that are located at a distance less than $\lambda > 0$.

Again, it should be noted that, i.e., when the PRCC is inviscid, when either the functions $\kappa_{i,j}$ are zero or the functions $g_k$ are constant, then, the only real-time measurement requirements of one vehicle are its own state and the positions of its adjacent vehicles (not their speeds, nor their orientations).

## 3.4. Statements of Main Results

The following Theorems guarantee that the closed-loop system (2.2) with (3.11), (3.12) and (2.2) with (3.22) and (3.23) satisfy properties (P1) and (P2).

**Theorem 1:** *For every $w_0 \in \Omega$ the initial-value problem (2.2) with (3.11), (3.12), (3.13) and $w(0) = w_0$ has a unique solution $w(t)$, defined for all $t \geq 0$, that satisfies $w(t) \in \Omega$ for all $t \geq 0$, as well as*

$$\lim_{t \to +\infty}(s_i(t)) = 0, \quad \lim_{t \to +\infty}\left(\frac{v_i(t)}{r_i(t)}\right) = \omega^*, \text{ for } i = 1,...,n \quad (3.28)$$

$$\lim_{t \to +\infty}(\dot{s}_i(t)) = 0, \quad \lim_{t \to +\infty}(F_i(t)) = 0, \text{ for } i = 1,...,n \quad (3.29)$$

*Furthermore, there exists a non-decreasing function $P_1 : \Re_+ \to \Re_+$ such that*

$$|F_i(t)| + |\delta_i(t)| \leq P_1(H(w_0)), \text{ for } t \geq 0, \ i = 1,...,n \quad (3.30)$$

**Theorem 2:** *For every $w_0 \in \Omega$ the initial-value problem (2.2) with (3.22), (3.23) and $w(0) = w_0$ has a unique solution $w(t)$, defined for all $t \geq 0$, that satisfies $w(t) \in \Omega$ for all $t \geq 0$, as well as*

$$\lim_{t \to +\infty}(s_i(t)) = 0, \quad \lim_{t \to +\infty}\left(\frac{v_i(t)}{r_i(t)}\right) = \omega^*, \text{ for } i = 1,...,n \quad (3.31)$$

$$\lim_{t \to +\infty}(\dot{s}_i(t)) = 0, \quad \lim_{t \to +\infty}(F_i(t)) = 0, \text{ for } i = 1,...,n \quad (3.32)$$

*Furthermore, there exists a non-decreasing function* $P_2 : \Re_+ \to \Re_+$ *such that*

$$|F_i(t)| + |\delta_i(t)| \leq P_2(H_R(w_0)), \text{ for } t \geq 0, \ i = 1,...,n \qquad (3.33)$$

**Remark:** (i) The results of Theorem 1 and Theorem 2 hold globally, i.e., for any initial condition $w_0 \in \Omega$.

(ii) The proofs of Theorem 1 and Theorem 2 rely on Barbălat's lemma ([19]) which only guarantees asymptotic convergence rate. Moreover, adding viscosity may increase the convergence rate of the CLFs. This is shown in the simulation scenarios of Section 5. For both the inviscid and viscous cases, the convergence rate is not exponential.

(iii) By selecting $R_{out} - R_{in}$ sufficiently small, the problem reduces to the study of vehicles moving on a circle as in traditional lane-based traffic.

(iv) Theorem 1 and Theorem 2 guarantee that property (P1) holds, namely, the following inequalities hold for all $i = 1,...,n$:

$$\begin{aligned}
R_{in} &< \inf_{t \geq 0}(r_i(t)) \leq \sup_{t \geq 0}(r_i(t)) < R_{out} \\
0 &< \inf_{t \geq 0}(v_i(t)) \leq \sup_{t \geq 0}(v_i(t)) < v_{\max} \\
-\Theta &< \inf_{t \geq 0}(s_i(t)) \leq \sup_{t \geq 0}(s_i(t)) < \Theta \\
&\inf_{t \geq 0}(d_{i,j}(t)) > L_{i,j}, \ i, j = 1,...,n, \ j \neq i
\end{aligned} \qquad (3.34)$$

## 4. Differences with the Straight-Road Case

Having designed cruise controllers for both the ring-road and the straight road (see [17], [18], [16]), it is possible to consider any other road represented by a closed or an open curve. Indeed, since a simple closed (open) curve is the image of a circle (line) under a homeomorphism (under a 1-1 continuous transformation with a continuous inverse), it is always possible to find a suitable change of coordinates to:
- Transform any non-self-intersecting curved road of constant width and infinite length into a straight road of constant width;
- Transform any road of constant width represented as a non-self-intersecting closed curve into a ring-road of constant width.

It should be noticed that it is not possible to find a homeomorphism to transform a straight road into a ring-road, and therefore, the cruise controllers (3.11), (3.12) and (3.22), (3.23) are not a direct consequence of the cruise controllers for the straight road under a change of coordinates. In addition, there are certain differences between the ring-road and the straight road listed next:

<u>D1.</u> *Geometry:* The straight road and the ring-road are *topologically* different. The straight road of constant width $2a > 0$ is a non-self-intersecting open curve defined on $\Re \times (-a, a)$, while the ring-road is a non-self-intersecting closed curve $\tilde{\gamma}$ defined on $[0, (R_{in} + R_{out})\pi] \times (R_{in}, R_{out})$ with $\tilde{\gamma}(0, \cdot) = \tilde{\gamma}((R_{in} + R_{out})\pi, \cdot)$. An immediate consequence of the geometry of the two cases is that the number of vehicles that operate on the ring-road is finite and depends on the width of the ring-road $R_{out} - R_{in}$, the length of the ring-road $(R_{in} + R_{out})\pi$, as well as the length of the vehicles $\sigma_i, i = 1,...,n$. On the contrary, the number of vehicles that can be placed in a straight road has no bound.

D2. *Speed set-point:* The objective of asymptotic convergence to a speed set-point in (P2) is different for each case due to the motion of vehicles of the ring-road and of the straight road. In particular, for the straight road case, the cruise controllers in [18] guarantee that the speed converges to a speed set-point $v^* \in (0, v_{max})$, i.e., $\lim_{t \to +\infty}(v_i(t)) = v^*$, whereas in the ring-road case, the NCC and PRCC controllers guarantee that the angular speed converges to an angular speed set-point $\omega^* \in \left(0, \frac{v_{max}}{R_{out}}\right)$, i.e., $\lim_{t \to +\infty}\left(\frac{v_i(t)}{r_i(t)}\right) = \omega^*$. While in both cases the speed remains bounded in $(0, v_{max})$, in the ring-road case we cannot regulate the speed $v_i$, but the angular speed $\frac{v_i}{r_i}$.

D3. *Steering angle limit:* The steering angle input $\delta_i$ satisfies $\lim_{t \to +\infty}\left(\delta_i(t) - \tan^{-1}\left(\frac{\sigma_i}{r_i(t)}\right)\right) = 0$ (recall (2.2) and (3.29), (3.32)) for the ring-road case contrary to the straight-road case where it holds that $\lim_{t \to +\infty}(\delta_i(t)) = 0$. More specifically, for the straight-road case, $\delta_i \equiv 0$ implies that the vehicle moves parallel to the road. For a vehicle operating on a ring-road however, its velocity vector and direction are always changing to follow a circular motion.

D4. *Distance between vehicles:* Since a straight road is defined on $\Re \times (-a, a)$, the inter-vehicle distance $d_{i,j} = \sqrt{(x_i - x_j)^2 + p_{i,j}(y_i - y_j)^2}$, $i = 1,...,n$, $j \neq i$, is unbounded, where $(x_i, y_i) \in \Re \times (-a, a)$ are the longitudinal and lateral position of vehicle $i$, respectively. However, for the case of a ring-road, the distance $d_{i,j}$ between vehicles defined by (2.5) is bounded, i.e., $d_{i,j} \leq 2R_{out}$ for all $i, j = 1,...,n$, $j \neq i$. Finally, although the selection of the constant $\lambda$ in (3.7) is equally important for both cases, since it defines the required real-time information zone around each vehicle (see discussion in [17]), the implications on the ring-road case are far greater. In particular, for selection of a very large $\lambda$, it is possible to receive information from all other vehicles in the ring-road. While this may increase the measurement requirements for the controllers, it may also lead to strong stabilizing effects of the density, see [15].

Besides the differences above, there are also certain similarities between the NCC and PRCC for the ring-road and the straight-road cases.

S1. *Nature of equilibria:* Consider system (2.2) under the change of coordinates $\tilde{\varphi}_i = \varphi_i - \omega^* t$, $t \geq 0$ where $\varphi_i$ is the angular coordinate. According to (2.2), (3.11), (3.12) and (2.2), (3.22), (3.23), the closed-loop systems are given by the same equations (2.2), (3.11), (3.12) and (2.2), (3.22), (3.23) with $\tilde{\varphi}_i$ replacing $\varphi_i$ everywhere (recall that $\varphi_i - \varphi_j = \tilde{\varphi}_i - \tilde{\varphi}_j$) and the differential equations $\dot{\tilde{\varphi}}_i = \frac{v_i}{r_i}\cos(s_i) - \omega^*$ in place of $\dot{\varphi}_i = \frac{v_i}{r_i}\cos(s_i)$ for $i = 1,...,n$. Dropping the tildes, we conclude that the closed-loop systems are given by the following equations for $i = 1,...,n$

$$\dot{r}_i = -v_i \sin(s_i)$$
$$\dot{\varphi}_i = \frac{v_i}{r_i}\cos(s_i) - \omega^*$$
$$\dot{s}_i = \frac{v_i}{\sigma_i}\tan(\delta_i) - \frac{v_i}{r_i}\cos(s_i)$$
$$\dot{v}_i = F_i$$
(4.1)

with (3.11), (3.12) or (3.22), (3.23). The state $w$ is still given by (2.6) and the state space $\Omega$ is given by (2.7).

Define the set:

$$E = \left\{ w \in \Omega : \begin{array}{l} v_i - \omega^* r_i = s_i = 0,\ i = 1,...,n \\ U'_i(r_i) + \sum_{j \neq i}\left(p_{i,j}(r_i - r_j) + r_j(1 - \cos(\varphi_i - \varphi_j))\right)d_{i,j}^{-1}V'_{i,j}(d_{i,j}) = 0,\ i = 1,...,n \\ \sum_{j \neq i} d_{i,j}^{-1}V'_{i,j}(d_{i,j})r_j \sin(\varphi_i - \varphi_j) = 0,\ i = 1,...,n \end{array} \right\} \quad (4.2)$$

The set $E$ defined by (4.2) is the set of equilibrium points of the closed-loop system with both the NCC and PRCC for the ring-road. In both cases this set is not bounded and not compact, since $\varphi_i$ are not bounded. Theorem 1 and Theorem 2 show that the solutions of both closed-loop systems tends to the set $E$ as $t \to +\infty$ (i.e., it can be shown that $\lim_{t \to +\infty}(dist(w(t),E)) = 0$; the set $E$ is a global attractor). However, Theorem 1 and Theorem 2 cannot guarantee that each solution tends to a single equilibrium point. Therefore, we cannot be sure that there exists an "ultimate" arrangement of the vehicles in the ring-road. The situation is utterly similar to the straight-road case (see [17], [18]). Advanced stability notions (e.g., stability with respect to two measures, see [32], [16]) are required in order to study the qualitative behavior of the solutions of both closed-loop systems.

S2. *Asymptotic convergence:* The NCC and PRCC (viscous and inviscid) for both the ring-road and the straight road only guarantee asymptotic convergence and not exponential convergence. The reason is that the proofs of Theorem 1 and Theorem 2 rely on Barbălat's Lemma, (see [19]) which only guarantees asymptotic convergence of the solution to the set $E$. In addition, the set of equilibrium points is an unbounded (and consequently, non-compact) set and may not present uniform attractivity properties, see [4].

S3. *State-Space*: In both cases the state-space $\Omega$ defined by (2.7) is an open set which may not be diffeomorphic to $\Re^{4n}$.

# 5. Numerical Simulations

In the simulation results below, we apply the proposed decentralized NCC and PRCC for both the viscous and inviscid case. Specifically, we consider a group of $n = 10$ vehicles on a lane-free ring-road with $R_{in} = 20$ and $R_{out} = 60$. The vehicle-repulsive potential functions $V_{i,j}$ and the boundary-repulsive potential function $U_i$ for both the NCC and PRCC are specified as

$$V_{i,j} = \begin{cases} q_1 \dfrac{(\lambda - d)^3}{d - L_{i,j}} & , L_{i,j} < d \leq \lambda \\ 0 & , d > \lambda \end{cases} \text{, for all } i, j = 1, \ldots, n \quad (5.1)$$

$$U_i(r) = \begin{cases} 0 & |r - R_m| \leq c \\ \dfrac{(r - R_m - c)^3 (r - R_m + c)^3}{(r - R_{in})(R_{out} - r)} & |r - R_m| > c \end{cases} \text{, for all } i = 1, \ldots, n \quad (5.2)$$

where $0 < c < \dfrac{R_{out} - R_{in}}{2}$, $q_1 > 0$ are design parameters, $R_m = \dfrac{R_{in} + R_{out}}{2}$, and $L_{i,j} = L$, $p_{i,j} = p$, $i, j = 1, \ldots, n$, $j \neq i$. Notice that $V_{i,j}$ and $U_i$ in (5.1) and (5.2) satisfy properties (3.1), (3.2), (3.3) and (3.4), (3.5), respectively. For small values of $q_1$, the values of $V_{i,j}$ (and consequently the acceleration $F_i$) will be smaller away from $L$, but will increase more sharply as $d$ approaches $L$. By adjusting the value of $c$, we can create an annulus $A = \{r \in (R_{in}, R_{out}) : R_m - c \leq r \leq R_m + c\}$ in the ring-road that satisfies $U_i(r) = 0$, $r \in A$, which may affect the final configuration of the vehicles relative to the boundaries $R_{in}$ and $R_{out}$. Notice that $V_{i,j}$ and $U_i$ above, satisfy (3.1), (3.2), (3.3), and (3.4), respectively. Finally, for both the NCC and PRCC we consider that the viscosity is given by

$$g_1(x) = g_2(x) = x$$
$$\kappa_{i,j}(d) = \begin{cases} q_2(\lambda - d)^2 & L < d \leq \lambda \\ 0 & d > \lambda \end{cases} \text{, for } i, j = 1, \ldots, n \quad (5.3)$$

where $q_2 \geq 0$ is a design parameter which can be selected to increase or decrease the effects of viscosity. Notice that if $q_2 \equiv 0$, then we obtain the inviscid NCC and PRCC cruise controllers. Moreover, we select for the NCC

$$f(x) = \dfrac{1}{2\varepsilon} \begin{cases} 0 & x \leq -\varepsilon \\ (x + \varepsilon)^2 & -\varepsilon < x < 0 \\ \varepsilon^2 + 2\varepsilon x & x \geq 0 \end{cases} \quad (5.4)$$

that satisfies (3.10), where $\varepsilon > 0$ is design parameter. Finally, for the PRCC we select $f_1(x) = \mu_1 x$ and $f_2(x) = \mu_2 x$, $x \in \Re$, where $\mu_1, \mu_2 > 0$.

To verify numerically and illustrate the results of Theorem 1, and Theorem 2, we assume that all vehicles have length $\sigma = 5$ and we set the angular speed set-point $\omega^* = 0.15$, the maximum speed $v_{max} = 10$ and select $\Theta = 0.17$ in order to satisfy condition (2.4) (recall that $R_{out} = 60$). Finally, we set $\varepsilon = 0.2$, $\mu_1 = 0.3$, $\mu_2 = 10^2$, $p_{i,j} = 5.11$, for all $i, j = 1, \ldots, n$, $q_1 = 3*10^{-3}$, $\lambda = 20$, $A = 0.5$, $b = 1$, $c = 10$, $q_2 = 0.1$ (for the viscous case) and $L = 6$.

Figure 2 shows the convergence of $\frac{v_i(t)}{r_i(t)}$ to $\omega^*$ for both the inviscid and viscous NCC. It is seen that by adding viscosity to the system, we have smoother convergence to $\omega^*$. This is also illustrated in Figure 3 which shows the evolution of the acceleration $\|F_i(t)\|_\infty$.

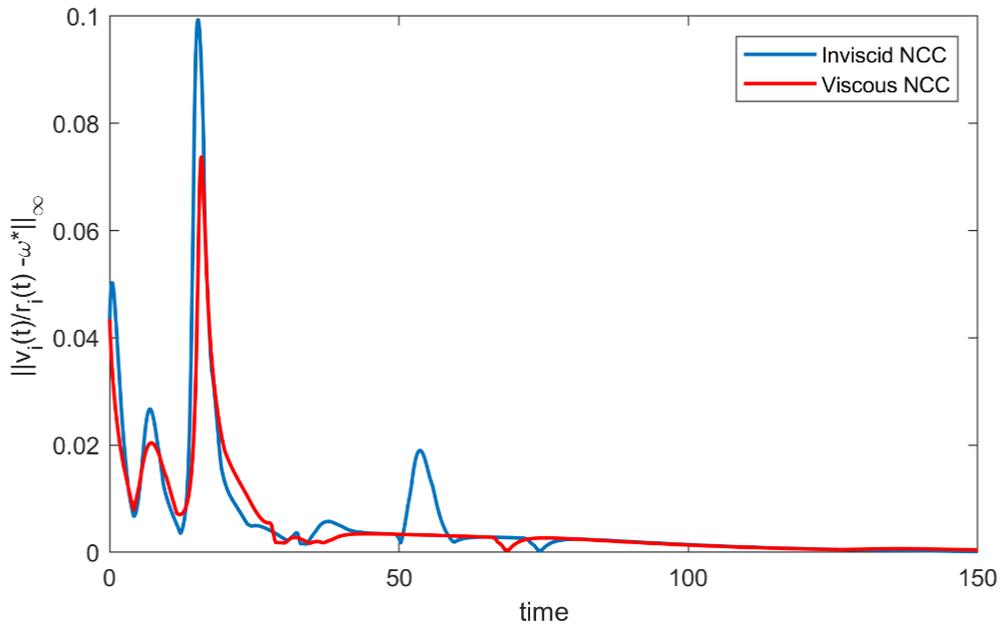

Figure 2: Evolution and convergence of $\left\|\frac{v_i(t)}{r_i(t)} - \omega^*\right\|_\infty$ for the Inviscid and Viscous NCC.

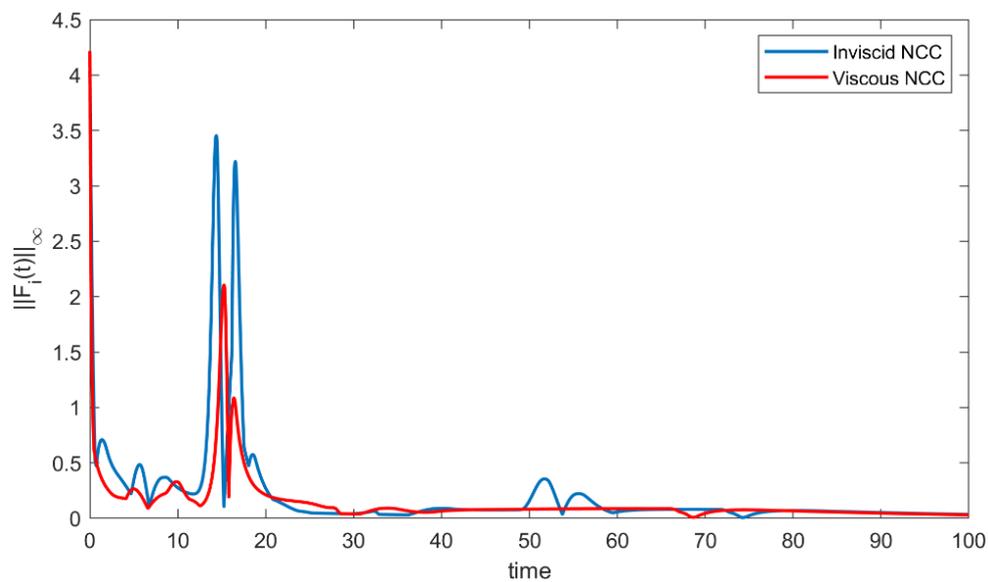

Figure 3: Evolution and convergence of $\|F_i(t)\|_\infty$ for the Inviscid and Viscous NCC.

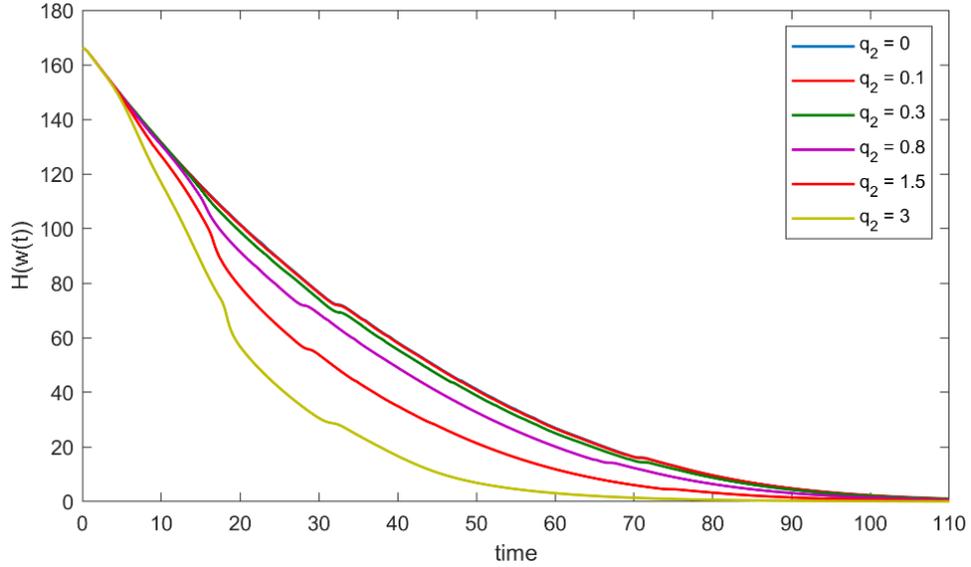

Figure 4: Evolution of the Control Lyapunov Function $H(w(t))$ for various $q_2$

Figure 4 shows the convergence of the Control Lyapunov Function $H(w)$ for the Inviscid and Viscous NCC. Finally, Figure 5 depicts the minimum inter-vehicle distance, verifying that there are no collisions among vehicles.

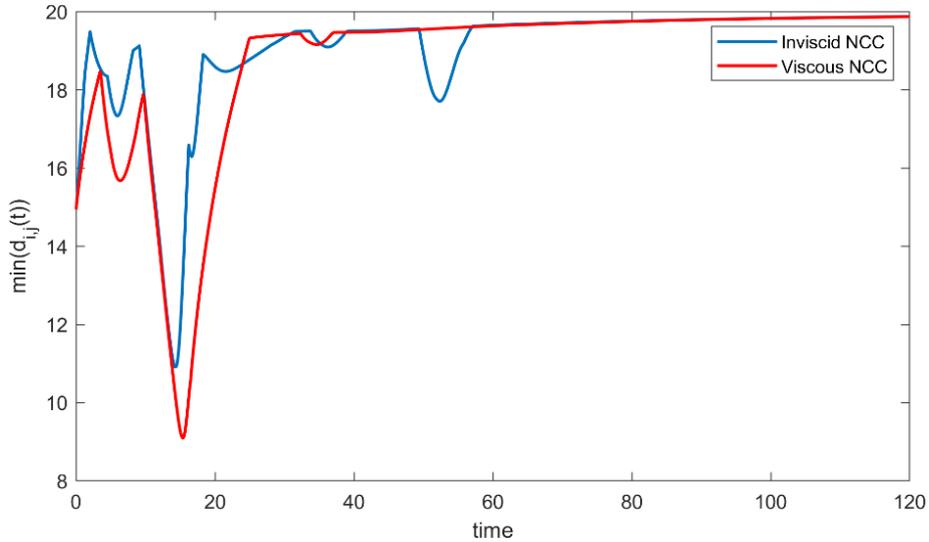

Figure 5: Minimum inter-vehicle distance for the inviscid and viscous NCC.

For the PRCC, we select $q_1 = 3 \cdot 10^{-5}$ and $q_2 = 0.1$ for the viscous case. Figure 6 shows the evolution and convergence of $\left\| \frac{v_i(t)}{r_i(t)} - \omega^* \right\|_\infty$ for the Inviscid and Viscous PRCC. Figure 7 shows the evolution of $\| F_i(t) \|_\infty$ for the viscous and inviscid PRCC. Figure 8 shows the minimum inter-vehicle distance, verifying that there are no collisions, and Figure 9 shows the convergence of $H_R(w(t))$.

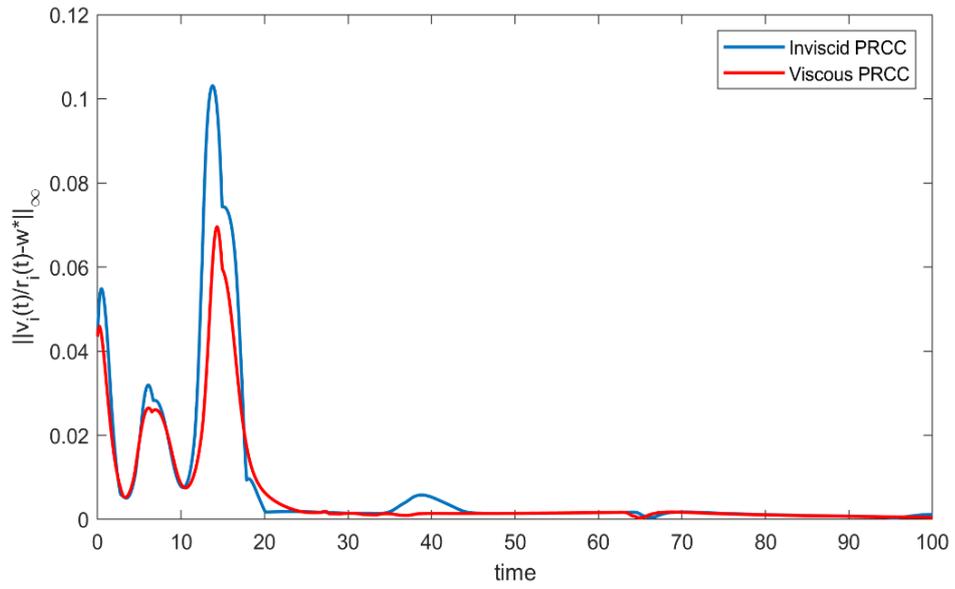

Figure 6: Evolution and convergence of $\left\|\dfrac{v_i(t)}{r_i(t)} - \omega^*\right\|_\infty$ for the Inviscid and Viscous NCC.

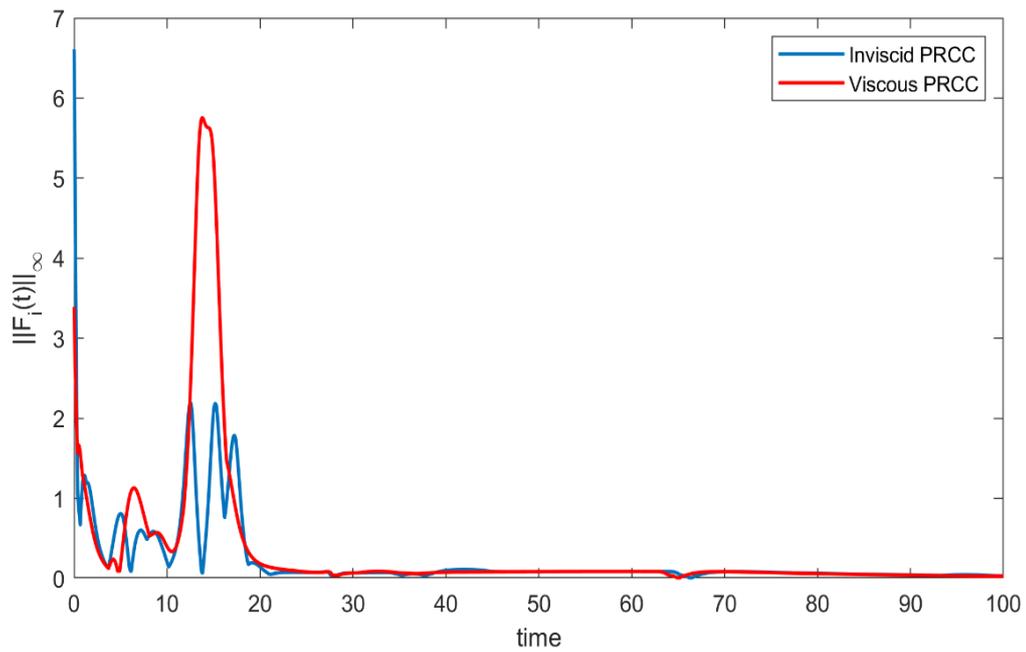

Figure 7: Evolution and convergence of $\left\|F_i(t)\right\|_\infty$ for the Inviscid and Viscous PRCC.

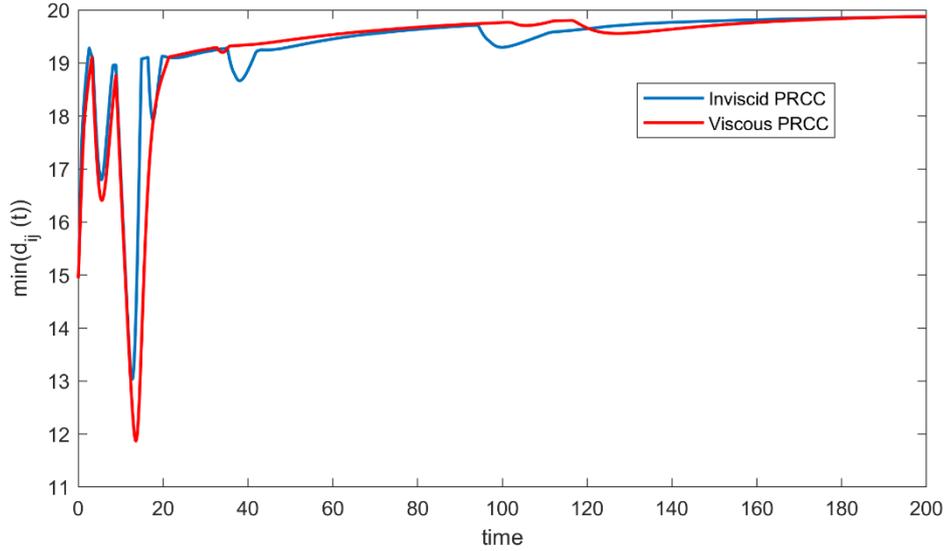

Figure 8: Minimum inter-vehicle distance for the inviscid and viscous PRCC

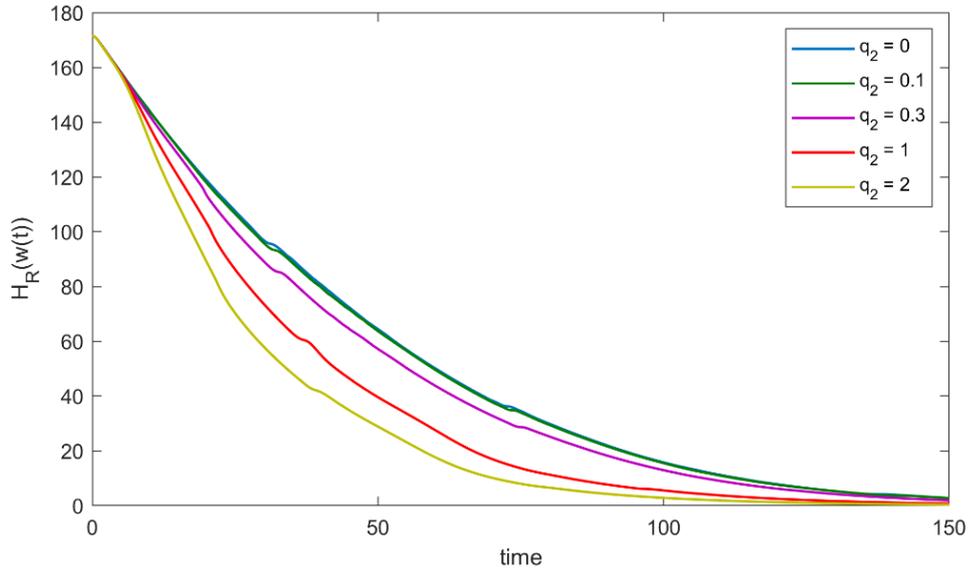

Figure 9: Evolution of Control Lyapunov Function $H_R(w(t))$ for various values of $q_2$.

## 6. Proofs

**Proof of Lemma 1:** For each $i = 1,...,n$, let $\xi_i \in \arg\min_{r \in (R_{in}, R_{out})}\{U_i(r)\}$, $u_i = U_i(\xi_i)$ and define

$$\alpha_i(s) := \min\{U_i(r) - u_i, r \in (R_{in}, s]\}, \quad s \in (R_{in}, \xi_i] \qquad (6.1)$$

$$\beta_i(s) := \min\{U_i(r) - u_i, r \in [\xi_i + s, R_{out})\}, \quad s \in [0, R_{out} - \xi_i) \qquad (6.2)$$

The family of functions $\alpha_i$, $i = 1,...,n$ are non-increasing with $a_i(\xi_i) = 0$, $\lim_{s \to R_{in}^+}(\alpha_i(s)) = +\infty$ and the functions $\beta_i$, $i = 1,...,n$ are non-decreasing with $\beta_i(0) = 0$ and $\lim_{s \to (R_{out} - \xi_i)^-}(\beta_i(s)) = +\infty$. We show first that for each $i = 1,...,n$ there exist constants $c_i^* \in (R_{in}, \xi_i]$, and decreasing functions $\zeta_i : \Re_+ \to (R_i, c_i^*]$ that satisfy $\zeta_i(H(w)) \le r_i$ in (3.9). Let $\bar{\alpha}_i(s) := \exp(R_{in} - s)\alpha_i(s)$, $s \in (R_{in}, \xi_i]$, which due to (6.1) satisfies $\lim_{s \to R_{in}^+} \bar{\alpha}_i(s) = +\infty$, $\bar{\alpha}_i(\xi_i) = 0$, and $\bar{\alpha}_i(s) \leq \alpha_i(s)$ for all $s > R_{in}$, $i = 1,...,n$. Let also $c_i^* = \inf \left\{ \arg\min_{r \in (R_{in}, \xi_i]} \{\bar{\alpha}_i(r)\} \right\}$ which satisfies $R_{in} < c_i^* \leq \xi_i$ since $\lim_{s \to R_{in}^+} \bar{\alpha}_i(s) = +\infty$. Moreover, $\bar{\alpha}_i(s)$ is decreasing on $(R_{in}, c_i^*]$. Consider the inverse functions $\zeta_i : \Re_+ \to (R_i, c_i^*]$, $i = 1,...,n$ of $\bar{\alpha}_i$ restricted on $(R_i, c_i^*]$ which are continuous and decreasing since $\bar{\alpha}_i$ are continuous and decreasing on $(R_i, c_i^*]$. Since definition (3.8) implies that $U_i(r_i) \leq H(w)$ for all $w \in \Omega$ and $i = 1,...,n$, it follows from definition (6.1) that $\bar{\alpha}_i(r_i) \leq \alpha_i(r_i) \leq H(w)$, $i = 1,...,n$. By taking into account that $\zeta_i$ are decreasing on $(R_{in}, c_i^*]$, the previous inequality shows that $r_i \geq \zeta_i(H(w)) > R_{in}$. We finally show that there exists a non-decreasing function $\eta_i$ such that $r_i \leq \eta_i(H(w))$ holds in (3.9). Lemma 1 in [17] suggests that there exist constants $s_i^* \in [0, R_{out} - \xi_i)$ and continuous, non-decreasing functions $\bar{\beta}_i : [0, R_{out} - \xi_i) \to \Re_+$, $i = 1,...,n$ with $\lim_{s \to (R_{out} - \xi_i)^-} \bar{\beta}_i(s) = +\infty$, $\bar{\beta}_i(s_i^*) = 0$, $\bar{\beta}_i(s) \leq \beta_i(s)$ for all $s \geq 0$ and $\bar{\beta}_i$ being increasing on $[s_i^*, R_{out} - \xi_i)$. Consider also the inverse functions $\eta_i : \Re_+ \to [s_i^* + \xi_i, R_{out})$, $i = 1,...,n$ of $\bar{\beta}_i(s - \xi_i)$ restricted on $[s_i^* + \xi_i, R_{out})$ which are continuous and increasing since $\bar{\beta}_i(s - \xi_i)$ are continuous and increasing. Moreover, notice that definition (3.8) implies that $U_i(r_i) \leq H(w)$ for all $w \in \Omega$ and $i = 1,...,n$. The previous inequality and the definition of $\eta_i$ above show that $r_i \leq \eta_i(H(w)) < R_{out}$. The existence of $\omega$ and $\rho_{i,j}$, $i, j = 1,...,n$, $j \neq i$ satisfying inequalities (3.9) is a minor modification of the proof of Proposition 1 in [17] and is omitted. The proof is complete. ◁

**Proof of Lemma 2:** The existence of the functions $\eta_i : \Re_+ \to [\xi_i, R_{out})$, $\zeta_i : \Re_+ \to (R_{in}, \xi_i]$, $\xi_i \in (R_{in}, R_{out})$, $\rho_{i,j} : \Re_+ \to (L_{i,j}, \lambda]$, $i, j = 1,...,n$, $i \neq j$, and $\omega : \Re_+ \to [0, \Theta)$ for which implication (3.20) holds, is a direct consequence of Lemma 1 with $H_R(w)$ in place of $H(w)$. Define

$$Q(v, r, s) := \frac{(v\cos(s) - r\omega^*)^2 + br^2v^2\sin(s)}{2r^2v(v_{\max} - v)}, \text{ for } v \in (0, v_{\max}), s \in [-\Theta, \Theta], r \in [R_{in}, R_{out}] \quad (6.3)$$

which satisfies $Q(v_i, r_i, s_i) \leq H_R(w)$ for all $w \in \Omega$ and $i = 1,...,n$ (a consequence of (3.19)) and $\omega^* \in \left(0, \frac{v_{\max}}{R_{out}}\right)$. Define for $v \in (0, v_{\max})$

$$f(v) := \min \{ Q(v, r, s) : s \in [-\Theta, \Theta], r \in [R_{in}, R_{out}] \} \quad (6.4)$$

which is well-defined due to compactness of $[-\Theta, \Theta] \times [R_{in}, R_{out}]$. Due to continuity of $Q(v, r, s)$ with respect to $v, r, s$ and compactness of $[-\Theta, \Theta]$ and $[R_{in}, R_{out}]$, it follows that $f$ defined in (6.4) is continuous, see Proposition 2.9 in [10].

We show next that $\lim_{v \to 0^+} f(v) = \lim_{v \to v_{\max}^-} f(v) = +\infty$. Consider an arbitrary sequence $\{v_i \in (0, v_{\max}) : i = 1,...\}$ with $\lim_{i \to +\infty}(v_i) = 0$. We next show that $\lim_{i \to +\infty}(f(v_i)) = +\infty$. Suppose on the contrary that there exists a bounded subsequence $\{f(v_k) \in (0, v_{\max}) : k = 1,...\}$, i.e., there exists $M > 0$ for which $f(v_k) \leq M$ (recall that $f(v) \geq 0$ for all $v \in (0, v_{max})$). Since the set $[-\Theta, \Theta] \times [R_{in}, R_{out}]$ is compact and due to continuity of $f$, there exists $(s_k, r_k) \in [-\Theta, \Theta] \times [R_{in}, R_{out}]$ for which $f(v_k) = Q(v_k, r_k, s_k)$. The latter, in conjunction with the facts that $f(v_k) = Q(v_k, r_k, s_k) \leq M$, $\lim_{k \to +\infty}(v_k) = 0$ and definition (6.3) implies that $\lim_{k \to +\infty}\left((v_k \cos(s_k) - r_k \omega^*)^2 + b v_k^2 r_k^2 \sin^2(s_k)\right) = 0$. By compactness of the set $[-\Theta, \Theta] \times [R_{in}, R_{out}]$ the sequence $\{(s_k, r_k) \in [-\Theta, \Theta] \times [R_{in}, R_{out}] : k = 1, 2, ...\}$ has a convergent subsequence $\{(s_m, r_m) \in [-\Theta, \Theta] \times [R_{in}, R_{out}] : m = 1, 2, ...\}$. In other words, there exist $(\bar{s}, \bar{r}) \in [-\Theta, \Theta] \times [R_{in}, R_{out}]$ with $\lim_{m \to +\infty}\left((s_m, r_m)\right) = (\bar{s}, \bar{r})$, $\lim_{m \to +\infty}(v_m) = 0$ and $\lim_{m \to +\infty}\left((v_m \cos(s_m) - r_m \omega^*)^2 + b v_m^2 r_m^2 \sin^2(s_m)\right) = (\bar{r} \omega^*)^2 = 0$. The latter equality implies that $\bar{r} = 0$; a contradiction with the fact that $\bar{r} \in [R_{in}, R_{out}]$. Thus, $\lim_{v \to 0^+}(f(v)) = +\infty$.

Using entirely analogous arguments as above and (2.4) it can be shown that $\lim_{v \to v_{max}^-}(f(v)) = +\infty$.

We proceed by showing the existence of $\chi, \psi$ that satisfy inequality (3.20). Notice first that for $\bar{v} \in \arg\min_{v \in (0, v_{max})} f(v)$ it holds that $f(\bar{v}) = 0$. Consider also the functions

$$p_1(s) = \min\{f(v), v \in (0, s]\}, \ s \in (0, \bar{v}] \quad (6.5)$$

$$p_2(s) = \min\{f(v), v \in [\bar{v} + s, v_{max})\}, \ s \in [0, v_{max} - \bar{v}) \quad (6.6)$$

The function $p_1$ is non-increasing with $p_1(\bar{v}) = 0$, $\lim_{s \to 0^+}(p_1(s)) = +\infty$ and the function $p_2$, is non-decreasing with $p_2(\bar{v}) = 0$ and $\lim_{s \to v_{max}^-}(p_2(s)) = +\infty$. Let $\bar{p}_1(s) := \exp(-s) p_1(s)$, $s \in (0, \bar{v}]$, which due to (6.5) satisfies $\lim_{s \to 0^+} \bar{p}_1(s) = +\infty$, $\bar{p}_1(\bar{v}) = 0$, and $\bar{p}_1(s) \leq p_1(s)$ for all $s > 0$. Moreover, $\bar{p}_1(s)$ is decreasing on $(0, \upsilon]$, $\upsilon = \inf\left\{\arg\min_{v \in (0, v_{max})}(f(v))\right\}$ with $\bar{p}_1(\upsilon) = 0$. Consider the inverse functions $\psi : \Re_+ \to (0, \upsilon]$, $i = 1, ..., n$ of $\bar{p}_1$ which is continuous and decreasing since $\bar{p}_1$ is continuous and decreasing. Since $\bar{p}_1(v_i) \leq p_1(v_i) \leq f(v_i) \leq H_R(w)$, if follows by definition of $\psi$, that $\psi(H_R(w)) \leq v_i$. Finally, by using Lemma 1 in [17], that there exists a constant $\bar{\upsilon} \in [0, v_{max} - \bar{v})$, and a continuous and non-decreasing function $\bar{p}_2 : [0, v_{max} - \bar{v}) \to \Re_+$, with $\lim_{s \to (v_{max} - \bar{v})^-} \bar{p}_2(s) = +\infty$, $\bar{p}_2(\bar{\upsilon}) = 0$, $\bar{p}_2(s) \leq p_2(s)$ for all $s \geq 0$ and $\bar{p}_2$ being increasing on $[\bar{\upsilon}, v_{max} - \bar{v})$. Consider also the inverse functions $\chi : \Re_+ \to [\bar{\upsilon}, v_{max})$ of $\bar{p}_2(s - \bar{\upsilon})$ restricted on $[\bar{\upsilon}, v_{max})$ which is continuous and increasing since $\bar{p}_2(s - \bar{\upsilon})$ is continuous and increasing. It follows by inequality $\bar{p}_2(v_i) \leq p_2(v_i) \leq f(v_i) \leq H_R(w)$ and definition of $\chi$ above, that $v_i \leq \chi(H_R(w))$ for all $w \in \Omega$. The proof is complete. ◁

The proof of Theorem 1 is performed by using Barbălat's lemma ([19]) and its following variant which uses uniform continuity of the derivative of a function.

**Lemma 3:** If a function $g \in C^2(\Re_+)$ satisfies $\lim_{t \to +\infty}(g(t)) \in \Re$ and $\sup_{t \geq 0}(|\ddot{g}(t)|) < +\infty$, then, $\lim_{t \to +\infty}(\dot{g}(t)) = 0$.

**Proof of Theorem 1:** We first show certain properties of $k_i(w)$ defined by (3.13) that will be used in the subsequent analysis. For every $w \in \Omega$, there exist non-decreasing functions $Q_i : \Re_+ \to \Re_+$, $i = 1,...,n$ such that the following inequalities hold

$$k_i(w) \geq \Phi(w) - G(w) \geq -\frac{k_i(w)}{r_i \omega^*}\left(v_{\max} \cos(s_i) - r_i \omega^*\right) \quad (6.7)$$

$$\mu_2 \leq k_i(w) \leq Q_i(w) \, , \, i = 1,...,n \quad (6.8)$$

where $\Phi$ and $G$ are defined by (3.16) and (3.17), respectively. Indeed, (6.7) is equivalent to the following inequality

$$\begin{aligned}k_i(w) &\geq \max\left\{\Phi(w) - G(w), -\frac{r_i \omega^*}{v_{\max} \cos(s_i) - r_i \omega^*}\left(\Phi(w) - G(w)\right)\right\} \\ &= \Phi(w) - G(w) + \max\left(0, \frac{v_{\max} \cos(s_i)}{v_{\max} \cos(s_i) - r_i \omega^*}\left(\Phi(w) - G(w)\right)\right)\end{aligned} \quad (6.9)$$

Inequality (6.7) is a direct consequence of (3.10), (3.13), and (6.9). Define next,

$$B_{i,j}(s) := \max\left\{|V'_{i,j}(d)| : s \leq d \leq \lambda\right\} \text{ for } s \in (L_{i,j}, \lambda] \, , \, j = 1,...,n, j \neq i \quad (6.10)$$

$$c_{i,j}(s) := \max\left\{\kappa_{i,j}(d) : s \leq d \leq \lambda\right\} \text{ for } s \in (L_{i,j}, \lambda], j = 1,...,n, j \neq i \quad (6.11)$$

$$\chi_i(s) := \max\left\{|U'_i(r)| : s \leq r \leq c_i^*\right\} \text{ for } s \in (R_{in}, c_i^*] \, , i = 1,...,n, \quad (6.12)$$

$$z_i(s) := \max\left\{|U'_i(r)| : c_i^* \leq r \leq s\right\} \text{ for } s \in [c_i^*, R_{out}) \, , i = 1,...,n, \quad (6.13)$$

where $c_i^* \in \arg\min_{r \in (R_{in}, R_{out})}\{|U'_i(r)|\}$, $i = 1,...,n$.

Moreover, for each $i = 1,...,n$, let $m_i \geq 2$ be the maximum number of points that can be placed within the area bounded by two concentric ellipses (in polar coordinates) with semi-major axes $\bar{L}_i = \min\{L_{i,j}, j = 1,...,n, j \neq i\}$ and $\lambda$ satisfying (3.7), and semi-minor axes $\frac{\bar{L}_i}{\max_{j \neq i} \sqrt{p_{i,j}}}$ and $\frac{\lambda}{\min_{j \neq i} \sqrt{p_{i,j}}}$ so that each point has distance (in the metric given by (2.5)) at least $\bar{L}_i$ from every other point. Then, it follows from (2.5), (3.2), the fact that $d_{i,j} > L_{i,j}$ for $i,j = 1,...,n$, $j \neq i$ and the definition of $m_i$ above, that the sums $\sum_{j \neq i} V'_{i,j}(d_{i,j}) \frac{r_j \sin(\varphi_i - \varphi_j)}{d_{i,j}}$, $\sum_{j \neq i}\left(p_{i,j}(r_i - r_j) + r_j(1 - \cos(\varphi_i - \varphi_j))\right)\frac{V'_{i,j}(d_{i,j})}{d_{i,j}}$ contain at most $m_i$ non-zero terms, namely the terms with $d_{i,j} \leq \lambda$. Definition (6.10), in conjunction with (3.9), and the fact that $p_{i,j} \geq 1$ for all $i,j = 1,...,n$, $j \neq i$, implies the following estimate holds for all $w \in \Omega$ and $i = 1,...,n$:

$$\max\left(\left|\frac{r_i}{\omega^*}\sum_{j\neq i}V'_{i,j}(d_{i,j})\frac{r_j \sin(\varphi_i - \varphi_j)}{d_{i,j}}\right|, \left|\sum_{j\neq i}\left(p_{i,j}(r_i - r_j) + r_j(1 - \cos(\varphi_i - \varphi_j))\right)\frac{V'_{i,j}(d_{i,j})}{d_{i,j}}\right|\right)$$

$$\leq \max\left(\frac{R_{out}^2}{\omega^*\bar{L}_i}\sum_{j\neq i}|V'(d_{i,j})|, \frac{2v_{max}}{\omega^*\bar{L}_i}\sum_{j\neq i}p_{i,j}|V'(d_{i,j})|\right)$$

$$\leq \frac{R_{out}^2 + 2v_{max}}{\omega^*\bar{L}_i}\max\left(\sum_{j\neq i}|V'(d_{i,j})|, \sum_{j\neq i}p_{i,j}|V'(d_{i,j})|\right) \quad (6.14)$$

$$\leq \frac{R_{out}^2 + 2v_{max}}{\omega^*\bar{L}_i}\sum_{j\neq i}p_{i,j}B_{i,j}(\rho_{i,j}(H(w))) \leq \frac{m_i(R_{out}^2 + 2v_{max})}{\omega^*\bar{L}_i}\max_{j\neq i}\left(p_{i,j}B_{i,j}(\rho_{i,j}(H(w)))\right)$$

Moreover, due to (3.6), (3.9) and (6.11), it holds that

$$\sum_{j\neq i}\kappa_{i,j}(d_{i,j}) \leq m_i \max_{j\neq i}\left(c_{i,j}(\rho_{i,j}(H(w)))\right). \quad (6.15)$$

It follows from definitions (3.16), (3.17), inequalities (6.14), (6.15), and the facts that $g_1$ is non-decreasing, $v_i \in (0, v_{max})$, $s_i \in (-\Theta, \Theta)$ and $r_i \in (R_{in}, R_{out})$, $i = 1,...,n$ that

$$|\Phi(w) - G(w)| \leq \frac{m_i(R_{out}^2 + 2v_{max})}{\omega^* \min_{j\neq i}(L_{i,j})}\max_{j\neq i}\left(p_{i,j}B_{i,j}(\rho_{i,j}(H(w)))\right)$$

$$+ \frac{m_i}{\omega^*}\max_{j\neq i}\left\{c_{i,j}(\rho_{i,j}(H(w)))\left(g_1\left(\frac{v_{max}}{R_{in}}\right) - g_1(0)\right)\right\} \quad (6.16)$$

Moreover, definition (3.8) implies that for all $w \in \Omega$ and for all $i = 1,...,n$ the following inequality holds (recall also that $\omega^* \in \left(0, \frac{v_{max}}{R_{out}}\right)$)

$$\frac{1}{v_{max}\cos(s_i) - r_i\omega^*} \leq \frac{A + (1-\cos(\Theta))H(w)}{A(v_{max} - r_i\omega^*) + (v_{max}\cos(\Theta) - r_i\omega^*)(1-\cos(\Theta))H(w)}$$

$$\leq \frac{A + (1-\cos(\Theta))H(w)}{A(v_{max} - R_{out}\omega^*) + (v_{max}\cos(\Theta) - R_{out}\omega^*)(1-\cos(\Theta))H(w)} \quad (6.17)$$

Inequality (6.8) is a direct consequence of (3.13), (6.16), and (6.17) with

$$Q_i(s) := \mu_2 + J_i(s) + \max\left(f(z) : |z| \leq \frac{v_{max}(A + (1-\cos(\Theta)))sJ_i(s)}{A(v_{max} - R_{out}\omega^*) + (v_{max}\cos(\Theta) - R_{out}\omega^*)(1-\cos(\Theta))s}\right) \quad (6.18)$$

and

$$J_i(s) := \frac{m_i\left(R_{out}^2 + 2v_{max}\right)}{\omega^* \min_{j \neq i}\left(L_{i,j}\right)} \max_{j \neq i}\left(p_{i,j}B_{i,j}\left(\rho_{i,j}(s)\right)\right) + \frac{m_i}{\omega^*} \max_{j \neq i}\left\{c_{i,j}(\rho_{i,j}(s))\left(g_1\left(\frac{v_{max}}{R_{in}}\right) - g_1(0)\right)\right\} \quad (6.19)$$

defined for $s \geq 0$, $i = 1, \ldots, n$.

We proceed now with the proof of Theorem 1. Let $w_0 \in \Omega$ and consider the unique solution $w(t)$ of the initial value problem (2.2), (3.11), (3.12) with initial condition $w(0) = w_0$. Using the fact that the set $\Omega$ is open (recall definitions (2.3), (2.7)), we conclude that there exists $t_{max} \in (0, +\infty]$ such that the solution $w(t)$ of (2.2), (3.11), (3.12) is defined on $[0, t_{max})$ and satisfies $w(t) \in \Omega$ for all $t \in [0, t_{max})$. Furthermore, if $t_{max} < +\infty$ then there exists an increasing sequence of times $\{t_i \in [0, t_{max}) : i = 1, 2, \ldots\}$ with $\lim_{i \to +\infty}(t_i) = t_{max}$ and either $\lim_{i \to +\infty}\left(dist\left(w(t_i), \partial\Omega\right)\right) = 0$ or $\lim_{i \to +\infty}\left(|w(t_i)|\right) = +\infty$.

Definitions (3.8), (3.14), (3.15), (3.16), and equations (2.2) give:

$$\nabla H(w)\dot{w} = \sum_{i=1}^{n} a(r_i, s_i, v_i)\dot{s}_i \sin(s_i) + \sum_{i=1}^{n}\left(b\dot{v}_i \sin(s_i) + \Lambda_i(w)\right)v_i \sin(s_i)$$
$$+ \sum_{i=1}^{n}\left(\frac{\dot{v}_i}{r_i}\cos(s_i) + \omega^*\Phi_i(w)\right)\left(\frac{v_i}{r_i}\cos(s_i) - \omega^*\right) \quad (6.20)$$
$$+ \omega^* \sum_{i=1}^{n}\sum_{j \neq i} V'_{i,j}(d_{i,j})\frac{r_i r_j \sin(\varphi_i - \varphi_j)}{d_{i,j}}$$

Due to (3.3) and (2.5) we have $\sum_{i=1}^{n}\sum_{j \neq i} V'_{i,j}(d_{i,j})\frac{r_i r_j \sin(\varphi_i - \varphi_j)}{d_{i,j}} = 0$. The latter, together with (3.11), (3.12), and (6.20) imply that

$$\nabla H(w)\dot{w} = -\mu_2 \sum_{i=1}^{n} \sin^2(s_i) - \sum_{i=1}^{n} k_i(w)\left(\frac{v_i}{r_i}\cos(s_i) - \omega^*\right)^2$$
$$- \frac{1}{2}\sum_{i=1}^{n}\sum_{j \neq i} \kappa_{i,j}(d_{i,j})\left(\sin(s_j) - \sin(s_i)\right)\left(g_2\left(\sin(s_j)\right) - g_2\left(\sin(s_i)\right)\right)$$
$$- \frac{1}{2}\sum_{i=1}^{n}\sum_{j \neq i} \kappa_{i,j}(d_{i,j})\left(\frac{v_j}{r_j}\cos(s_j) - \frac{v_i}{r_i}\cos(s_i)\right)\left(g_1\left(\frac{v_j}{r_j}\cos(s_j)\right) - g_1\left(\frac{v_i}{r_i}\cos(s_i)\right)\right)$$
$$\leq -\mu_2 \sum_{i=1}^{n} \sin^2(s_i) - \mu_1 \sum_{i=1}^{n}\left(\frac{v_i}{r_i}\cos(s_i) - \omega^*\right)^2$$

$$\text{for all } w \in \Omega \quad (6.21)$$

It follows from (6.21) that for all $w \in \Omega$

$$\nabla H(w)\dot{w} \leq 0 \quad (6.22)$$

Since $w(t) \in \Omega$ for all $t \in [0, t_{max})$, it follows from (2.3) and (2.7) that $v_i(t) \in (0, v_{max})$, $s_i(t) \in (-\Theta, \Theta)$ and $r_i(t) \in (R_{in}, R_{out})$ for all $t \in [0, t_{max})$ and $i = 1,...,n$. Thus, (2.2) implies that $0 \leq \dot{\varphi}_i(t) \leq \frac{v_{max}}{R_{in}}$ for all $t \in [0, v_{max})$ and $i = 1,...,n$. In addition, inequality (6.22) implies that

$$H(w(t)) \leq H(w_0) \text{ for all } t \in [0, t_{max}) \tag{6.23}$$

Consequently, we obtain from (3.9) and (6.23) that for all $t \in [0, t_{max})$ and $i, j = 1,...,n$, $j \neq i$

$$|s_i(t)| \leq \omega(H(w_0)) < \Theta, R_{in} < \zeta_i(H(w_0)) \leq r_i(t) \leq \eta_i(H(w_0)) < R_{out}$$
$$d_{i,j}(t) \geq \rho_{i,j}(H(w_0)) > L_{i,j}, \varphi_i(0) \leq \varphi_i(t) \leq \varphi_i(0) + \frac{v_{max}}{R_{in}} t \tag{6.24}$$

From (6.8), (6.23) and (6.24) we further get that

$$\mu_2 < k_i(w(t)) \leq Q_i(H(w_0)) \text{ for all } t \in [0, t_{max}), i = 1,...,n \tag{6.25}$$

where $Q_i$ are defined in (6.18). Inequality (6.7), definition (3.11) and (2.2) imply that the following inequalities hold for all $t \in [0, v_{max})$ and $i = 1,...,n$

$$k_i(w(t))(v_{max} - v_i(t)) \geq \dot{v}_i(t) \geq -k_i(w(t))v_i(t) \tag{6.26}$$

Using (6.24) and (6.26) we have for all $t \in [0, t_{max})$ and $i = 1,...,n$, the following inequalities

$$v_i(0) \exp(-Q_i(H(w_0))t) + (1 - \exp(-Q_i(H(w_0))t)v_{max}$$
$$\geq v_i(t) \geq v_i(0) \exp(-Q_i(H(w_0))t) \tag{6.27}$$

Suppose that $t_{max} < +\infty$. Inequalities (6.23), (6.24), (6.27) and definitions (2.3) and (2.7) imply that for every increasing sequence of times $\{t_i \in [0, t_{max}) : i = 1, 2, ...\}$ with $\lim_{i \to +\infty}(t_i) = t_{max}$ we cannot have $\lim_{i \to +\infty}(dist(w(t_i), \partial\Omega)) = 0$ or $\lim_{i \to +\infty}(|w(t_i)|) = +\infty$. Thus, $t_{max} = +\infty$.

We show next that (3.30) holds. Since $w(t) \in \Omega$ for all $t \geq 0$, due to (2.3), (6.8), (6.23) and (6.26) it follows that

$$|\dot{v}_i(t)| \leq F_i(t) \leq Q_i(H(w_0))v_{max} \tag{6.28}$$

where $Q_i$ are defined in (6.18). Inequality (3.30) is a direct consequence of definition (3.12) and inequality (6.28) with $P_1(s) = \max_{i=1,...,n} \left\{ Q_i(s)v_{max} + \frac{\pi}{2} \right\}$, $s \geq 0$.

In order to show (3.28), we have that definition (3.14), (3.15), (3.18), inequalities (6.12), (6.13), (6.14), (6.15), (6.23), (6.24), and the facts that, $\omega^* \in \left(0, \frac{v_{max}}{R_{out}}\right)$, $g_2$ is non-decreasing, and $\sin(s)$ is increasing for $s \in \left(-\frac{\pi}{2}, \frac{\pi}{2}\right)$, imply that for all $w \in \Omega$ and $i, j = 1,...,n$, $j \neq i$, the following inequalities hold

$$|\Lambda_i(w)| \leq \left(\frac{v_{max}}{R_{in}} - \omega^*\right)\frac{v_{max}}{R_{in}^2} + \max\left(\chi_i(\zeta_i(H(w_0))), z_i(\eta_i(H(w_0)))\right)$$
$$+ \frac{m_i\left(R_{out}^2 + 2v_{max}\right)}{\omega^* \bar{L}_i} \max_{j \neq i}\left(p_{i,j} B_{i,j}\left(\rho_{i,j}(H(w_0))\right)\right) \tag{6.29}$$

$$|M_i(w)| \leq m_i \max_{j \neq i}\left(c_{i,j}(\rho_{i,j}(H(w_0)))\right)(g_2(\sin(\Theta)) - g_2(\sin(-\Theta))) \tag{6.30}$$

$$a(r,s,v) \geq \frac{A}{(1-\cos(\Theta))^2} \tag{6.31}$$

Inequalities (6.24), (6.28), (6.29), (6.30), (6.31) and definition (3.12) imply that $\tan(\delta_i(t))$ is bounded for all $t \geq 0$ and $i = 1,...,n$. The latter, together with (2.2) and the fact that $v_i(t) \in (0, v_{max})$, $r_i(t) \in (R_{in}, R_{out})$ and $s_i(t) \in (-\Theta, \Theta)$, $i = 1,...,n$ imply that $\dot{r}_i(t)$ and $\dot{s}_i(t)$ are bounded for all $i = 1,...,n$. Thus, using (2.2) (3.30), and (6.28) we conclude that there exists a constant $\xi > 0$ such that

$$\frac{d}{dt}\left(\mu_2 \sum_{i=1}^n \sin^2(s_i(t)) + \mu_1 \sum_{i=1}^n \left(\frac{v_i(t)}{r_i(t)}\cos(s_i(t)) - \omega^*\right)^2\right) \leq \xi \tag{6.32}$$

Finally, (6.21) implies that

$$\int_0^\infty \mu_2 \sum_{i=1}^n \sin^2(s_i(t)) + \mu_1 \sum_{i=1}^n \left(\frac{v_i(t)}{r_i(t)}\cos(s_i(t)) - \omega^*\right)^2 dt \leq H(w_0) \tag{6.33}$$

It follows from (6.32), (6.33) and Barbălat's lemma ([19]) that (3.28) holds.

Finally, we show that (3.29) holds for the solutions $w(t)$ of (2.2), (3.11), (3.12) by exploiting Lemma 3 with $g(t) = v_i(t)$ and $g(t) = s_i(t)$ for $i = 1,...,n$. Since (3.28) holds, it suffices to show that $\ddot{s}_i(t)$ and $\dot{F}_i(t)$ are bounded.

We show first that $\dot{F}_i(t)$ is bounded for all $i = 1,...,n$. Notice that, since $v_i(t) \in (0, v_{max})$, $s_i(t) \in (-\Theta, \Theta)$, $r_i(t) \in (R_{in}, R_{out})$ for all $t \geq 0$, $i = 1,...,n$ we have, due to continuity of $g'_k$, $k = 1, 2$, that

$$\left|g'_1\left(\frac{v_i(t)}{r_i(t)}\cos(s_i(t))\right)\right| \leq \max_{x \in \left[0, \frac{v_{max}}{R_{in}}\right]}(g'_1(x)), \quad \text{for } i = 1,...,n \text{ and } t \geq 0 \tag{6.34}$$

$$\left|g'_2(\sin(s_i(t)))\right| \leq \max_{x \in [-1,1]}(g'_2(x)), \quad \text{for } i = 1,...,n \text{ and } t \geq 0. \tag{6.35}$$

Define also the non-increasing functions

$$\bar{\kappa}_{i,j}(s) := \max\left\{|\kappa'_{i,j}(d)| : s \leq d \leq \lambda\right\} \text{ for } s \in (L_{i,j}, \lambda], i, j = 1,...,n, j \neq i. \tag{6.36}$$

where $\kappa'_{i,j}$ is the derivative of the function $\kappa_{i,j}$ satisfying (3.5), (3.6). Then, from (3.6), (6.24), and (6.36) we get

$$\left|\kappa'_{i,j}(d_{i,j}(t))\right| \le \overline{\kappa}_{i,j}(\rho_{i,j}(H(w_0))) \; , \; i,j=1,\ldots,n, \; j \ne i. \qquad (6.37)$$

Notice also that due to (6.24), boundedness of $\dot{r}_i(t)$, $\dot{\varphi}_i(t)$ for $i=1,\ldots,n$ (recall (2.2), (6.24) and the fact that $v_i(t) \in (0, v_{max})$) and definition (2.5) it follows that $\dot{d}_{i,j}(t)$ is bounded. Moreover, notice that since (3.2) and (6.24) hold, it follows that $V'_{i,j}(d_{i,j}(t)), V''_{i,j}(d_{i,j}(t))$ are bounded for all $i,j=1,\ldots,n$ with $j \ne i$. Due to (6.24), boundedness of $V'_{i,j}(d_{i,j}(t)), V''_{i,j}(d_{i,j}(t))$, $\dot{d}_{i,j}(t)$, definition (3.16), (3.17) and inequalities (6.37) and (6.14) it follows that $\frac{d}{dt}(\Phi_i(w(t)) - G_i(w(t)))$ is bounded for all $i=1,\ldots,n$.

We show next that $\frac{d}{dt}(k_i(w(t)))$ is bounded. We have from definition (3.13),

$$\frac{d}{dt}(k_i(w(t))) = \frac{d}{dt}(\Phi_i(w(t)) - G_i(w(t))) + f'\left(-\frac{v_{max}\cos(s_i(t))}{v_{max}\cos(s_i(t)) - r_i(t)\omega^*}(\Phi_i(w(t)) - G_i(w(t)))\right) \\ \times \frac{d}{dt}\left(-\frac{v_{max}\cos(s_i(t))}{v_{max}\cos(s_i(t)) - r_i(t)\omega^*}(\Phi_i(w(t)) - G_i(w(t)))\right) \qquad (6.38)$$

Formula (6.38), together with (6.24), (6.16), (6.19), boundedness of $\frac{d}{dt}(\Phi_i(w(t)) - G_i(w(t)))$ and inequality

$$\left|f'\left(-\frac{v_{max}\cos(s_i(t))}{v_{max}\cos(s_i(t)) - r_i(t)\omega^*}(\Phi_i(w(t)) - G_i(w(t)))\right)\right| \\ \le \max\left(|f'(z)| : |z| \le \frac{v_{max}(A + (1-\cos(\Theta)))H(w_0)J(H(w_0))}{A(v_{max} - R_{out}\omega^*) + (v_{max}\cos(\Theta) - R_{out}\omega^*)(1-\cos(\Theta))H(w_0)}\right)$$

imply that each $\frac{d}{dt}(k_i(w(t)))$, $i=1,\ldots,n$ is bounded.

Finally, using boundedness of $\dot{r}_i(t)$, $\dot{s}_i(t)$, $\frac{d}{dt}(\Phi_i(w(t)) - G_i(w(t)))$, $\frac{d}{dt}(k_i(w(t)))$, inequalities (6.16), (6.24), (6.25), (6.28) and formula

$$\dot{F}_i(t) = -\frac{d}{dt}(k_i(w(t)))\left(v_i(t) - \frac{r_i(t)\omega^*}{\cos(s_i(t))}\right) - k_i(w(t))\frac{d}{dt}\left(v_i(t) - \frac{r_i(t)\omega^*}{\cos(s_i(t))}\right) \\ -\frac{d}{dt}\left(\frac{r_i(t)\omega^*}{\cos(s_i(t))}\right)(\Phi_i(w(t)) - G_i(w(t))) - \frac{r_i(t)\omega^*}{\cos(s_i(t))}\frac{d}{dt}(\Phi_i(w(t)) - G_i(w(t)))$$

we conclude that $\dot{F}_i(t)$, $i=1,\ldots,n$ are bounded.

We show next that $\ddot{s}_i(t)$ is bounded for each $i=1,\ldots,n$. Due to (2.2), boundedness of $\dot{r}_i(t)$, $\dot{v}_i(t)$, $\dot{s}_i(t)$, $\tan(\delta_i(t))$, $i=1,\ldots,n$, definition (3.12), the facts that $v_i(t) \in (0, v_{max})$, $r_i(t) \in (R_{in}, R_{out})$, $s_i(t) \in (-\Theta, \Theta)$, $i=1,\ldots,n$, inequalities (6.24), and formula

$$\ddot{s}_i(t) = \frac{\dot{v}_i(t)}{\sigma_i}\tan(\delta_i(t)) + \frac{v_i(t)}{\sigma_i}\frac{d}{dt}\big(\tan(\delta_i(t))\big) - \frac{d}{dt}\left(\frac{v_i(t)}{r_i(t)}\cos(s_i(t))\right) \qquad (6.39)$$

it suffices to show that

$$\begin{aligned}\frac{d}{dt}(\tan(\delta_i(t))) &= \frac{d}{dt}\left(\frac{\sigma_i}{r_i(t)}\cos(s_i(t))\right) \\ &\quad - \frac{\sigma_i}{a(r_i(t),s_i(t),v_i(t))v_i(t)}\Big(\mu_2\sin(s_i(t)) - \big(bF_i(t)\sin(s_i(t)) + \Lambda_i(w(t))\big)v_i(t) - M_i(w(t))\Big)\end{aligned} \qquad (6.40)$$

is bounded for all $i=1,\ldots,n$.

Boundedness of $\frac{d}{dt}\left(\frac{\sigma_i}{r_i(t)}\cos(s_i(t))\right)$, $\frac{d}{dt}\left(\frac{\sigma_i}{a(r_i(t),s_i(t),v_i(t))v_i(t)}\right)$, $\frac{d}{dt}(\sin(s_i(t)))$, is a direct consequence of (2.2), (3.14), (6.24), (6.28), boundedness of $\dot{r}_i(t)$, $\dot{s}_i(t)$ for all $i=1,\ldots,n$. Moreover, due to (3.2) and (6.24), we have that $U_i'(r_i(t))$, $U_i''(r_i(t))$, $V_{i,j}'(d_{i,j}(t))$, and $V_{i,j}''(d_{i,j}(t))$ are bounded for all $i,j=1,\ldots,n$, $j\neq i$. Thus, definitions (3.15), (3.16), (3.18) boundedness of $\dot{d}_{i,j}(t)$, and inequalities (6.24), (6.35), and (6.37) imply that $\frac{d}{dt}(\Lambda_i(w(t)))$ and $\frac{d}{dt}(M_i(w(t)))$ are bounded. Therefore, according to (6.40), boundedness of $\dot{F}_i(t)$, and inequalities (6.29), (6.30), and (6.31), we have that $\ddot{s}_i(t)$ is bounded for all $i=1,\ldots,n$ as well. This completes the proof. ◁

**Proof of Theorem 2:** Let $w_0 \in \Omega$ and consider the unique solution $w(t)$ of the initial value problem (2.2), (3.22), (3.23) with $w(0)=w_0$. Using the fact that the set $\Omega$ is open (recall definitions (2.3), (2.7)), we conclude that there exists $t_{\max} \in (0,+\infty]$ such that the solution $w(t)$ of (2.2), (3.22), (3.23) is defined on $[0,t_{\max})$ and satisfies $w(t) \in \Omega$ for all $t \in [0,t_{\max})$. Furthermore, if $t_{\max}<+\infty$ then there exists an increasing sequence of times $\{t_i \in [0,t_{\max}): i=1,2,\ldots\}$ with $\lim_{i\to+\infty}(t_i)=t_{\max}$ and either $\lim_{i\to+\infty}\big(dist(w(t_i),\partial\Omega)\big)=0$ or $\lim_{i\to+\infty}(|w(t_i)|)=+\infty$.

By using (2.2), (2.5), and (3.19), definition (3.16), (3.24), (3.25), (3.26), and (3.27), it follows that for all $w \in \Omega$

$$\begin{aligned}\nabla H_R(w)\dot{w} &= \sum_{i=1}^n \gamma(r_i,s_i,v_i)\dot{s}_i\sin(s_i) + \sum_{i=1}^n \big(\dot{v}_i\zeta(s_i,v_i)+Z_i(w)\big)v_i\sin(s_i) \\ &\quad + \sum_{i=1}^n \big(\dot{v}_i q(r_i,s_i,v_i)+\omega^*\Phi_i(w)\big)\left(\frac{v_i}{r_i}\cos(s_i)-\omega^*\right) \\ &\quad + \omega^* \sum_{i=1}^n\sum_{j\neq i} V_{i,j}'(d_{i,j})\frac{r_i r_j \sin(\varphi_i-\varphi_j)}{d_{i,j}}\end{aligned} \qquad (6.41)$$

Moreover, due to (2.5) and (3.3) we have $\sum_{i=1}^{n}\sum_{j\neq i}V'_{i,j}(d_{i,j})\frac{r_i r_j \sin(\varphi_i - \varphi_j)}{d_{i,j}} = 0$. The latter, together with (6.41), (3.17), (3.18), (3.22), (3.23), and the fact that $g_k$, $k=1,2$ are non-decreasing, imply that for all $w \in \Omega$ the following inequality holds

$$\nabla H_R(w)\dot{w} = -\sum_{i=1}^{n}\left(\frac{v_i}{r_i}\cos(s_i) - \omega^*\right)f_1\left(\frac{v_i}{r_i}\cos(s_i) - \omega^*\right) - \sum_{i=1}^{n}\sin(s_i)f_2(\sin(s_i))$$
$$-\frac{1}{2}\sum_{i=1}^{n}\sum_{j\neq i}\kappa_{i,j}(d_{i,j})\left(\sin(s_j) - \sin(s_i)\right)\left(g_2(\sin(s_j)) - g_2(\sin(s_i))\right)$$
$$-\frac{1}{2}\sum_{i=1}^{n}\sum_{j\neq i}\kappa_{i,j}(d_{i,j})\left(\frac{v_j}{r_j}\cos(s_j) - \frac{v_i}{r_i}\cos(s_i)\right)\left(g_1\left(\frac{v_j}{r_j}\cos(s_j)\right) - g_1\left(\frac{v_i}{r_i}\cos(s_i)\right)\right)$$
$$\leq -\sum_{i=1}^{n}\left(\frac{v_i}{r_i}\cos(s_i) - \omega^*\right)f_1\left(\frac{v_i}{r_i}\cos(s_i) - \omega^*\right) - \sum_{i=1}^{n}\sin(s_i)f_2(\sin(s_i))$$

(6.42)

It follows from (3.21), and (6.42) that

$$\nabla H_R(w)\dot{w} \leq 0, \text{ for all } w \in \Omega \tag{6.43}$$

Since $w(t) \in \Omega$ for all $t \in [0, t_{max})$, it follows from (2.3), (2.7) that $v_i(t) \in (0, v_{max})$, $s_i(t) \in (-\Theta, \Theta)$, $r_i(t) \in (R_{in}, R_{out})$ for all $t \in [0, t_{max})$ and $i = 1, ..., n$. Thus, (2.2) implies that $0 \leq \dot{\varphi}_i(t) \leq \frac{v_{max}}{R_{out}}$ for all $t \in [0, t_{max})$ and $i = 1, ..., n$. Moreover, inequality (6.43) implies that

$$H_R(w(t)) \leq H_R(w_0), \text{ for all } t \in [0, t_{max}) \tag{6.44}$$

Consequently, we obtain from (6.44) and (3.20) that for all $t \in [0, t_{max})$ and $i, j = 1, ..., n$, $j \neq i$:

$$\eta_i(H_R(w_0)) \leq r_i(t) \leq \zeta_i(H_R(w_0)),$$
$$\ell(H_R(w_0)) \leq v_i(t) \leq \psi(H_R(w_0)), |s_i(t)| \leq \omega(H_R(w_0)),$$
$$\varphi_i(0) \leq \varphi_i(t) \leq \varphi_i(0) + \frac{v_{max}}{R_{out}}t, d_{i,j}(t) \geq \rho_{i,j}(H_R(w_0)) > L_{i,j}$$

(6.45)

which imply that for every increasing sequence of times $\{t_i \in [0, t_{max}) : i = 1, 2, ...\}$ with $\lim_{i \to +\infty}(t_i) = t_{max}$ we cannot have $\lim_{i \to +\infty}(dist(w(t_i), \partial\Omega)) = 0$ or $\lim_{i \to +\infty}(|w_i(t)|) = +\infty$. Therefore, $t_{max} = \infty$.

We show next that (3.31) and (3.32) hold. Similar to (6.14), and by using (6.45) we have that

$$\max\left(\left|\frac{r_i}{\omega^*}\sum_{j\neq i}V'_{i,j}(d_{i,j})\frac{r_j\sin(\varphi_i - \varphi_j)}{d_{i,j}}\right|, \left|\sum_{j\neq i}\left(p_{i,j}(r_i - r_j) + r_j(1 - \cos(\varphi_i - \varphi_j))\right)\frac{V'_{i,j}(d_{i,j})}{d_{i,j}}\right|\right)$$
$$\leq \frac{R_{out}^2 + 2v_{max}}{\omega^*\overline{L}_i}\sum_{j\neq i}p_{i,j}B_{i,j}(\rho_{i,j}(H_R(w))) \leq \frac{m_i(R_{out}^2 + 2v_{max})}{\omega^*\overline{L}_i}\max_{j\neq i}\left(p_{i,j}B_{i,j}(\rho_{i,j}(H_R(w)))\right)$$

(6.46)

where $\bar{L}_i = \min_{j \neq i}(L_{i,j})$. Moreover, due to (3.6), (6.11), and (6.45), it holds that

$$\sum_{j \neq i} \kappa_{i,j}(d_{i,j}) \leq m_i \max_{j \neq i}\left(c_{i,j}(\rho_{i,j}(H_R(w)))\right). \tag{6.47}$$

Notice also that due to (2.4) and the facts that $v_i \in (0, v_{max})$, $r_i \in (R_{in}, R_{out})$, $\cos(s_i) > \cos(\Theta)$, $s_i \in (-\Theta, \Theta)$, $\Theta \in \left(0, \frac{\pi}{2}\right]$, $i = 1,...,n$ we have that

$$v_{max} v_i \cos(s_i) - 2 r_i v_i \omega^* + r_i \omega^* v_{max} > R_{in} \omega^* (v_{max} - v_i) \tag{6.48}$$

Inequality (6.48), definition (3.24) and inclusions $v \in (0, v_{max})$, $s \in (-\Theta, \Theta)$, $r \in (R_{in}, R_{out})$ imply that

$$\frac{1}{q(r,s,v)} \leq \frac{2 R_{out} v_{max}^3}{R_{in} \omega^*} \tag{6.49}$$

Moreover, from (3.25), inequality $b > \dfrac{1}{R_{in}^2}$ (recall 3.19) and (6.45), we also obtain the following inequality

$$\frac{1}{\gamma(r,s,v)} \leq \frac{(1-\cos(\Theta))^2}{A} \tag{6.50}$$

Using (3.26), (6.45) and the facts that $\ell$, is non-decreasing and $\psi$ is non-increasing, we obtain the following estimate,

$$\zeta(s_i(t), v_i(t)) \leq \frac{b v_{max}}{2(v_{max} - \psi(H_R(w_0)))^2 \, \ell(H_R(w_0))}. \tag{6.51}$$

Moreover, due to the facts that $v_i(t) \in (0, v_{max})$, $r_i(t) \in (R_{in}, R_{out})$, $s_i(t) \in \left(-\frac{\pi}{2}, \frac{\pi}{2}\right)$ for $i = 1,...,n$ and continuity of $f_k$, $k = 1, 2$, we have that

$$\left| f_1\left(\frac{v_i(t)}{r_i(t)} \cos(s_i(t)) - \omega^*\right) \right| \leq \xi_1 := \max_{x \in \left[-\omega^*, \frac{v_{max}}{R_{in}} - \omega^*\right]} (|f_1(x)|), \text{ for } i = 1,...,n \text{ and } t \geq 0 \tag{6.52}$$

$$|f_2(\sin(s_i(t)))| \leq \xi_2 := \max_{x \in [-1,1]} (|f_2(x)|), \text{ for } i = 1,...,n \text{ and } t \geq 0. \tag{6.53}$$

Then, from (3.22), (6.49), (6.52), and with analogous arguments as in (6.16), we obtain the following estimate for $i = 1,...,n$ and $t \geq 0$

$$|F_i(t)| \le \frac{2R_{out}v_{max}^3}{R_{in}\omega^*}\left(\xi_1 + \frac{m_i\left(R_{out}^2 + 2v_{max}\right)}{\min_{j\ne i}(L_{i,j})}\max_{j\ne i}\left(p_{i,j}B_{i,j}\left(\rho_{i,j}(H_R(w))\right)\right)\right.$$
$$\left.+m_i\max_{j\ne i}\left\{c_{i,j}(\rho_{i,j}(H_R(w)))\left(g_1\left(\frac{v_{max}}{R_{in}}\right) - g_1(0)\right)\right\}\right) \quad (6.54)$$

Inequality (3.33) is a direct consequence of (6.54) and definition (3.23) with

$$P_2(s) = \max_{i=1,\ldots,n}\left(\frac{2R_{out}v_{max}^3}{R_{in}\omega^*}\left(\xi_1 + \frac{m_i\left(R_{out}^2 + 2v_{max}\right)}{\min_{j\ne i}(L_{i,j})}\max_{j\ne i}\left(p_{i,j}B_{i,j}\left(\rho_{i,j}(s)\right)\right)\right.\right.$$
$$\left.\left.+m_i\max_{j\ne i}\left\{c_{i,j}(\rho_{i,j}(s))\left(g_1\left(\frac{v_{max}}{R_{in}}\right) - g_1(0)\right)\right\}\right) + \frac{\pi}{2}\right), \quad s \ge 0$$

From the facts that $v_i \in (0, v_{max})$, $r_i \in (R_{in}, R_{out})$, $s_i \in (-\Theta, \Theta)$, $i=1,\ldots,n$, inequalities (6.44), (6.47), (6.54), (3.18), (3.20), (3.27), and definitions (6.12), (6.13) we have the following estimates

$$|Z_i(w)| \le \left(\frac{v_{max}}{R_{in}} - \omega^*\right)\frac{1}{(v_{max} - \psi(H_R(w_0)))R_{in}^2} + \max\left(\chi_i(\zeta_i(H_R(w_0))), z_i(\eta_i(H_R(w_0)))\right)$$
$$+ \frac{m_i\left(R_{out}^2 + 2v_{max}\right)}{\omega^*\overline{L}_i}\max_{j\ne i}\left(p_{i,j}B_{i,j}\left(\rho_{i,j}(H_R(w_0))\right)\right) \quad (6.55)$$

$$|M_i(w)| \le m_i\max_{j\ne i}\left(c_{i,j}(\rho_{i,j}(H_R(w_0)))\right)(g_2(\sin(\Theta)) - g_2(\sin(-\Theta))) \quad (6.56)$$

From inequalities (6.45), (6.54), (6.55), (6.56), (6.50), (6.51), (6.53) and definition (3.23) we obtain that $\tan(\delta_i(t))$ is bounded for all $i=1,\ldots,n$. The latter, in conjunction with (2.2) and the fact that $v_i(t) \in (0, v_{max})$, $r_i \in (R_{in}, R_{out})$, $s_i \in (-\Theta, \Theta)$, $i=1,\ldots,n$ implies that $\dot{r}_i(t)$, $\dot{\varphi}_i(t)$, and $\dot{s}_i(t)$ are bounded for all $i=1,\ldots,n$. Moreover, since $f_j \in C^1(\Re)$, $j=1,2$ and due to inclusions $v_i(t) \in (0, v_{max})$, $s_i(t) \in (-\Theta, \Theta)$, and $r_i(t) \in (R_{in}, R_{out})$, $i=1,\ldots,n$ we also obtain that

$$\left|f_1'\left(\frac{v_i(t)}{r_i(t)}\cos(s_i(t)) - \omega^*\right)\right| \le \max_{x\in\left[-\omega^*, \frac{v_{max}}{R_{in}} - \omega^*\right]}\left(|f_1'(x)|\right), \text{ for } i=1,\ldots,n \text{ and } t \ge 0 \quad (6.57)$$

$$\left|f_2'(\sin(s_i(t)))\right| \le \max_{x\in[-1,1]}\left(|f_2'(x)|\right), \text{ for } i=1,\ldots,n \text{ and } t \ge 0. \quad (6.58)$$

Thus, using (2.2), (6.45), (6.54), (6.57), (6.58), and boundedness of $\dot{r}_i(t)$ and $\dot{s}_i(t)$, we conclude that there exists a constant $\xi > 0$ such that

$$\frac{d}{dt}\left(\sum_{i=1}^n\left(\frac{v_i(t)}{r_i(t)}\cos(s_i(t)) - \omega^*\right)f_1\left(\frac{v_i(t)}{r_i(t)}\cos(s_i(t)) - \omega^*\right) + \sum_{i=1}^n\sin(s_i(t))f_2(\sin(s_i(t)))\right) \le \xi \quad (6.59)$$

Finally, (6.42) implies that

$$\int_0^\infty \sum_{i=1}^n \left(\frac{v_i(t)}{r_i(t)}\cos(s_i(t)) - \omega^*\right) f_1\left(\frac{v_i(t)}{r_i(t)}\cos(s_i(t)) - \omega^*\right) + \sum_{i=1}^n \sin(s_i(t)) f_2(\sin(s_i(t))) dt \le H(w_0) \quad (6.60)$$

It follows from (6.59), (6.60) and Barbălat's lemma ([19]) that (3.31) holds.

Finally, we show that (3.29) holds for the solutions $w(t)$ of (2.2), (3.22), (3.23) by exploiting Lemma 2 with $g(t) = v_i(t)$ and $g(t) = s_i(t)$ for $i = 1,...,n$. Since (3.31) holds, it suffices to show that $\ddot{s}_i(t)$ and $\dot{F}_i(t)$ are bounded.

We show first that $\dot{F}_i(t)$ is bounded for all $i = 1,...,n$. Notice that due to (6.45), boundedness of $\dot{r}_i(t)$, $\dot{\varphi}_i(t)$ for $i = 1,...,n$ and definition (2.5) it follows that $\dot{d}_{i,j}(t)$ is bounded. Moreover, notice that since (3.2) and (6.45) hold, it follows that $V'_{i,j}(d_{i,j}(t)), V''_{i,j}(d_{i,j}(t))$ are bounded for all $i, j = 1,...,n$ with $j \ne i$. Furthermore, from (3.6), (6.44), (6.45), and (6.36) we have that

$$\left|\kappa'(d_{i,j}(t))\right| \le \bar{\kappa}_{i,j}(\rho_{i,j}(H_R(w_0))), \quad i, j = 1,...,n, j \ne i. \quad (6.61)$$

Due to (6.45), boundedness of $V'_{i,j}(d_{i,j}(t)), V''_{i,j}(d_{i,j}(t))$, $\dot{d}_{i,j}(t)$, definitions (3.16), (3.17) and inequalities (6.46) and (6.61) it follows that $\frac{d}{dt}(\Phi_i(w(t)) - G_i(w(t)))$ is bounded for all $i = 1,...,n$.

Next, due to (6.45), boundedness of $\dot{s}_i(t)$, $\dot{r}_i(t)$, (6.54), and definition (3.24) it follows that $\frac{d}{dt}\left(\frac{1}{q(r_i(t), s_i(t), v_i(t))}\right)$ is bounded for each $i = 1,...,n$. Finally, boundedness of $\dot{F}_i(t)$ is a direct consequence of (3.22), (6.45), (6.49), (6.52), (6.57), and boundedness of $\Phi_i(w) - G_i(w)$, $\frac{d}{dt}\left(\frac{1}{q(r_i(t), s_i(t), v_i(t))}\right)$, and $\frac{d}{dt}(\Phi_i(w(t)) - G_i(w(t)))$ for $i = 1,...,n$.

We show next that $\ddot{s}_i(t)$ is bounded for each $i = 1,...,n$. Due to (2.2), boundedness of $\dot{r}_i(t)$, $\dot{v}_i(t)$, $\dot{s}_i(t)$, $\tan(\delta_i(t))$, $i = 1,...,n$, definition (3.23), and formula (6.39) it suffices to show that

$$\frac{d}{dt}\left(\frac{\sigma_i}{r_i}\cos(s_i) - \frac{\sigma_i}{\gamma(r_i, s_i, v_i)v_i}\left(f_2(\sin(s_i)) - (\zeta_i(s_i, v_i)F_i + Z_i(w))v_i - M_i(w)\right)\right) \quad (6.62)$$

is bounded for all $i = 1,...,n$.

First notice that boundedness of $\frac{d}{dt}\left(\frac{\sigma_i}{r_i(t)}\cos(s_i(t))\right)$, $\frac{d}{dt}\left(\frac{\sigma_i}{\gamma(r_i(t), s_i(t), v_i(t))v_i(t)}\right)$, $\frac{d}{dt}(f_2(\sin(s_i(t))))$, $\frac{d}{dt}(\zeta_i(s_i(t), v_i(t)))$ is a direct consequence of (2.2), (3.24), (3.25), (3.26), (6.45), (6.54), (6.58), and of the facts that $v_i(t) \in (0, v_{\max})$, $s_i(t) \in (-\Theta, \Theta)$, and $r_i(t) \in (R_{in}, R_{out})$ for all $i = 1,...,n$. Moreover, due to (3.2) and (6.45), we have that $U'_i(r_i(t))$, $U''_i(r_i(t))$, $V'_{i,j}(d_{i,j}(t))$, and $V''_{i,j}(d_{i,j}(t))$ are bounded for all $i, j = 1,...,n$, $j \ne i$. Thus, definitions (3.18), (3.27), boundedness of $\dot{d}_{i,j}(t)$, and inequalities (6.45), (6.35), and (6.37) imply that $\frac{d}{dt}(Z_i(w(t)))$ and $\frac{d}{dt}(M_i(w(t)))$ are bounded. Finally, boundedness of all the previous terms implies that (6.62) is bounded for all $i = 1,...,n$. This completes the proof. ◁

# 7. Conclusions

The paper introduced two families of cruise controllers for autonomous vehicles operating on lane-free ring-roads. By expressing the Control Lyapunov Functions on measures of the energy of the system with the kinetic energy expressed in ways similar to Newtonian or relativistic mechanics, we derived decentralized feedback laws (cruise controllers) that guarantee collision avoidance between vehicles and with the boundary of the ring-road; that the speeds of all vehicles are always positive and remain below a given speed limit; and that all vehicle angular speeds converge to a given speed set-point.

The differences and similarities of the cruise controllers for the straight road and the ring-road are also discussed. Having designed controllers for both lane-free ring-roads and straight roads, we can consider any road represented by an open or closed non-self-intersecting curve of constant width by appropriate change of coordinates.

## Acknowledgments

The research leading to these results has received funding from the European Research Council under the European Union's Horizon 2020 Research and Innovation programme/ ERC Grant Agreement n. [833915], project TrafficFluid.